\documentclass[10pt, leqno]{article}
\usepackage{amssymb, amsbsy, amsmath, amsfonts, amssymb, amscd, mathrsfs}

\usepackage[english]{babel}
\usepackage[T1]{fontenc}
\usepackage{indentfirst}

\usepackage{ifpdf}
\ifpdf
\usepackage[pdftex]{color, graphicx}
\else
\usepackage[dvips]{color, graphicx} 
\fi

\usepackage{caption}

\makeatletter
\@addtoreset{equation}{section}
\makeatother

\newtheorem{statement}{}[section]
\newtheorem{theoreme}[statement]{Theorem}
\newtheorem{lemme}[statement]{Lemma}

\newtheorem{proposition}[statement]{Proposition}

\newtheorem{corollaire}[statement]{Corollary}

\newcommand\C{\mathbb C}
\newcommand\N{\mathbb N}
\newcommand\R{\mathbb R}
\newcommand\T{\mathbb T}
\newcommand\D{\mathbb D}

\newcommand\e{{\rm e}}
\renewcommand\P{\mathbb P}
\newcommand\esp{\mathbb E}
\newcommand\eps{\varepsilon}
\newcommand\ind{{\rm 1\kern-.30em I}}
\newcommand\qed{\hfill $\square$}
\renewcommand \Re{{\mathfrak R}{\rm e}\,}
\renewcommand \Im{{\mathfrak I}{\rm m}\,}

\let\phi=\varphi

\title{Compact composition operators on Bergman-Orlicz spaces}
\author{\it Pascal Lef\`evre, Daniel Li,\\ \it Herv\'e Queff\'elec, Luis Rodr{\'\i}guez-Piazza}

\date{\footnotesize \today}

\begin{document}

\maketitle

\noindent{\bf Abstract.} \emph{
We construct an analytic self-map $\phi$ of the unit disk and an Orlicz 
function $\Psi$ for which the composition operator of symbol $\phi$ is compact on the Hardy-Orlicz space 
$H^\Psi$, but not on the Bergman-Orlicz space ${\mathfrak B}^\Psi$. For that, we first 
prove a Carleson embedding theorem, and then characterize the compactness of composition operators on 
Bergman-Orlicz spaces, in terms of Carleson function (of order $2$). We show that this Carleson 
function is equivalent to the Nevanlinna counting function of order $2$.}
\medskip

\medskip



\section{Introduction and notation}

\subsection{Introduction}

Due to the Littlewood subordination principle, the boundedness of composition operators $C_\phi$, defined by 
$C_\phi (f) = f \circ \phi$, on Hardy spaces $H^p$, as well as on Bergman spaces ${\mathfrak B}^p$, 
$1 \leq p \leq \infty$, is automatic. Their compactness is something much more subtle, but is well 
understood now, and there are two well-separated cases. First, the case $p = \infty$, for which 
$C_\phi \colon H^\infty \to H^\infty$ is compact if and only if $\| \phi \|_\infty < 1$ (note that 
${\mathfrak B}^\infty = H^\infty$). Secondly, the case $p < \infty$, for which the compactness does 
not depend on $p$. For Hardy spaces, this fact, proved by J. Shapiro and P. Taylor (\cite{Shapiro-Taylor}), is 
not completely trivial, and is due to the good factorization properties of functions in $H^p$. For the scale of 
Bergman spaces ${\mathfrak B}^p$, the factorization properties are not so good, but the independence with 
respect to $p$ follows from the following characterization (\cite{McCluer-Shapiro}, Corollary~4.4):  for 
$1 \leq p < \infty$, $C_\phi \colon {\mathfrak B}^p \to {\mathfrak B}^p$ is compact if and only if the pull-back 
measure of the area-measure ${\cal A}$ by $\phi$ is a vanishing $2$-Carleson measure. The case $p = 2$ 
(proved in \cite{Boyd}) gives then, for $1 \leq p < \infty$:
\begin{equation}\label{compacite sur $B^p$}
C_\phi \colon {\mathfrak B}^p \to {\mathfrak B}^p \text{ is compact } \quad \Longleftrightarrow \quad 
\lim_{|z| \to 1} \frac {1 - |\phi (z)|} {1 - |z|} = \infty\,.
\end{equation}
\par

In both cases (Hardy and Bergman), a brutal change of situation occurs when we pass from finite values 
of $p$ to the value $p = \infty$, and the need was felt for an intermediate scale between $H^p$ and 
$H^\infty$, or between ${\mathfrak B}^p$ and ${\mathfrak B}^\infty$. This is what we did (\cite{CompOrli}), 
in full detail, with Hardy-Orlicz spaces $H^\Psi$ associated with an Orlicz function $\Psi$ (and began to do 
for Bergman-Orlicz spaces ${\mathfrak B}^\Psi$). We introduced a generalization of the notion of Carleson 
measure, and proved a contractivity property of those Carleson measures $m_\phi = \phi^\ast (m)$, 
attached to an analytic self-map $\phi \colon \D \to \D$, which turned out to be central to obtain a 
necessary and sufficient condition for the compactness of $C_\phi \colon H^\Psi \to H^\Psi$. In that paper, 
we also began a study of the compactness of composition operators on ${\mathfrak B}^\Psi$. We proved, 
in particular, but implicitly (see the comments at the beginning of Section~\ref {Hardy-Bergman}), 
that, if the Orlicz function $\Psi$ grows very fast (satisfying the so-called $\Delta^2$ condition), 
then the compactness of $C_\phi \colon H^\Psi \to H^\Psi$ implies its compactness as an operator 
$C_\phi \colon {\mathfrak B}^\Psi \to {\mathfrak B}^\Psi$. On the other hand, it is well-known that the 
compactness on $H^p$ implies the compactness on ${\mathfrak B}^p$ because it is easy to see that the right-hand 
side of \eqref{compacite sur $B^p$} is implied by the compactness on $H^p$. One might think that 
it is generally easier to achieve compactness on ${\mathfrak B}^\Psi$ than on $H^\Psi$. The main result of the 
present work is the existence of an analytic self-map $\phi$ of $\D$ and an Orlicz function $\Psi$ such that 
the composition operator $C_\phi$ is compact on $H^\Psi$ but not on ${\mathfrak B}^\Psi$. For that, we first 
have to characterize the compactness of composition operators on Bergman-Orlicz spaces. More 
precisely, the paper is organized as follows.\par
\medskip

In Section~\ref{Carleson embeddings}, given two Orlicz functions $\Psi_1$ and $\Psi_2$, and a finite positive 
measure $\mu$ on the unit disk $\D$, we investigate under which conditions the canonical inclusion 
$I_\mu \colon {\mathfrak B}^{\Psi_1} \to L^{\Psi_2} (\mu)$, defined by $I_\mu (f) = f$, is either bounded, 
or compact. In Theorem~\ref{Theorem boundedness}, we give a \emph{necessary} condition and a 
\emph{sufficient} condition, in terms of the Carleson function $\rho_\mu$ of $\mu$, 
for the boundedness of $I_\mu$. Analoguously, we have a similar statement 
(Theorem~\ref{Theorem compactness}) for the compactness of $I_\mu$. In general, these necessary and those 
sufficient conditions do not fit.\par

In Section~\ref{compactness composition op}, we prove one of the main results of this paper 
(Theorem~\ref{homogeneite}) under the form of a contractivity principle for the pull-back measure 
${\cal A}_\phi$ of the planar Lebesgue measure ${\cal A}$ on $\D$ by $\phi$. The proof is rather long and 
uses a Calder{\'o}n-Zygmund decomposition, as well as an elementary, but very useful, inequality due to 
Paley and Zygmund. This contractivity principle eliminates the absence of fitness mentioned above and 
allows us to have a \emph{necessary and sufficient} condition for the compactness of 
$C_\phi \colon {\mathfrak B}^\Psi \to {\mathfrak B}^\Psi$ in terms of the same Carleson function 
$\rho_{{\cal A}_\phi} = \rho_{\phi, 2}$ (Theorem~\ref{CNS compactness composition operators}).\par

In Subsection~\ref{Nevanlinna}, we consider the Nevanlinna counting function $N_{\phi, 2}$ (initiated in 
\cite{Shapiro}), adapted to the Bergman case, and we compare it with the $2$-Carleson 
function $\rho_{\phi, 2}$ of $\phi$. These two functions turn out to be equivalent, in the sense precised 
in Theorem~\ref{equiv N-C Bergman}. This extends to the Bergman case (and follows from) such an equivalence 
for the Hardy case, that we recently established in \cite{LLQR-N}, Theorem~1.1.\par

Finally, in Section~\ref{Hardy-Bergman}, we exploit the necessary and sufficient conditions that we 
established, either on $H^\Psi$ and on ${\mathfrak B}^\Psi$, to give (Theorem~\ref{theo Hardy-Bergman}) 
an example of an analytic self-map $\phi \colon \D \to \D$ and of a fairly irregularly varying Orlicz function 
$\Psi$ such that, contrary to the general intuition, $C_\phi \colon H^\Psi \to H^\Psi$ is compact, whereas 
$C_\phi \colon {\mathfrak B}^\Psi \to {\mathfrak B}^\Psi$ is not compact. This is due to the fact that we can 
evaluate, in an accurate way, the two Carleson functions $\rho_\phi$ and $\rho_{\phi, 2}$ of $\phi$.

\bigskip
\noindent{\bf Acknowledgement.} The fourth-named author is partially supported by a Spanish research 
project MTM2006-05622.


\subsection{Notation}

We shall denote by $\D$ the open unit disk $\{z \in \C\,; \ |z| < 1\}$ of the complex plane, and its 
boundary, the unit circle, by $\T$. The normalized area measure $d {\cal A} = dx\, dy/\pi$ on $\D$ will be 
denoted by ${\cal A}$.\par
\smallskip

For any $\xi \in \T$, we define, for $0 < h < 1$, the Carleson window $W (\xi, h)$ by 
\begin{displaymath}
W (\xi, h) = \{ z\in \D\,; \ |z| \geq 1 - h \quad \text{and} \quad | \arg (z \overline{\xi} ) | \leq \pi h \}. 
\end{displaymath}
We shall also use the ``circular'' Carleson windows $S (\xi, h)$ defined by 
$S (\xi, h) = \{ z\in \D\,; \ |z - \xi | < h\}$. Since $S (\xi, h) \subseteq W (\xi, h) \subseteq S (\xi, 5h)$, the 
measures of $W (\xi, h)$ and of $S (\xi, h)$ are equivalent, up to constants.
\par\smallskip

For any finite positive measure $\mu$ on $\D$, we define, for $0 < h \leq 1$, the \emph{Carleson function} 
of $\mu$ by:
\begin{equation}
\rho_\mu (h) = \sup_{|\xi| = 1} \mu \big( W (\xi, h) \big) \,,
\end{equation}
and we set:
\begin{equation}
K_{\mu, 2} (h) = \sup_{0 < t < h} \frac {\rho_\mu (t)} {t^2}
\end{equation}

When $\rho_\mu (h) = O\,(h^2)$, one says that $\mu$ is a \emph{$2$-Carleson measure}; we also say that 
$\mu$ is a \emph{Bergman-Carleson measure}, to insist that the order $2$ is adapted to the Bergman spaces. 
When $\mu = {\cal A}_\phi$ is the pull-back measure of ${\cal A}$ by an analytic self-map 
$\phi \colon \D \to \D$, we shall simply write $\rho_{\phi, 2}$ and $K_{\phi, 2}$ instead of 
$\rho_{{\cal A}_\phi}$ and $K_{{\cal A}_\phi, 2}$ respectively. We shall say that $\rho_{\phi, 2}$ the 
\emph{$2$-dimensional Carleson function} of $\phi$.\par
\smallskip

The Hastings-Luecking sets of size $2^{-n}$ are defined by:
\begin{displaymath}
\Delta_k =  \Big\{z\in \D\,;\ 1 - \frac {1} {2^n} \leq |z| < 1 - \frac {1} {2^{n + 1} } \ \text{and}\  
\frac{(2j - 1) \pi}{2^n} \leq \arg z < \frac{(2j+ 1) \pi}{2^n}\,\Big\} \,,
\end{displaymath}
where $k = 2^n + j - 1$, $n \geq 0$, $0 \leq j \leq 2^n - 1$ (note that $\Delta_0 = D (0, 1/2)$).\par
\medskip

An Orlicz function $\Psi$ is a positive increasing convex function $\Psi \colon [0, \infty) \to [0, \infty)$ 
such that $\Psi (0) = 0$ (and $\Psi (\infty) = \infty$). If $\mu$ is a positive measure on $\D$, the Orlicz 
space $L^\Psi (\mu)$ is the space of (classes of) measurable functions $f \colon \D \to \C$ such that 
$\int_\D \Psi (|f|/ C)\, d{\cal A} < \infty$, for some constant $ C > 0$, and the norm $\| f \|_\Psi$ is defined 
as the infimum of all constants $C > 0$ for which $\int_\D \Psi (|f|/ C)\, d{\cal A} \leq 1$. The 
\emph{Bergman-Orlicz space} is the subspace of  $L^\Psi ({\cal A})$ whose members are analytic in $\D$.\par
The Hardy-Orlicz space $H^\Psi$ is the subspace of  $H^1$ whose boundary values are in the Orlicz 
space $L^\Psi (\T, m)$.
\bigskip
  
We refer to \cite{Duren-S} (see also \cite{HKZ}, and \cite{Zhu}) for the theory of Bergman spaces and to 
\cite{Rao} for more information about Orlicz spaces.


\section{Carleson embeddings}\label{Carleson embeddings}

We consider in this Section the ``embedding'' map 
$I_\mu \colon {\mathfrak B}^{\Psi_1} \to L^{\Psi_2} (\mu)$, defined by $I_\mu (f) = f$, where  $\mu$ is an 
arbitrary finite positive Borel measure on $\D$ and $\Psi_1$ and $\Psi_2$ are two Orlicz functions.


\subsection{Boundedness}\label{boundedness}

\begin{theoreme}\label{Theorem boundedness}
Given $\mu$ a finite positive Borel measure on $\D$ and $\Psi_1$ and $\Psi_2$ two Orlicz functions, let 
$I_\mu \colon {\mathfrak B}^{\Psi_1} \to L^{\Psi_2} (\mu)$ be the canonical map defined by 
$I_\mu (f) = f$. One has:\par

1) If $I_\mu$ is bounded, then there is a constant $A >0$ such 
that: 
\begin{equation}\label{CN continue}
\qquad \qquad \rho_\mu (h) \leq \frac {1} {\Psi_2 [A \Psi_1^{-1} (1/h^2) ] } 
\, \raise 1,5pt \hbox{,} \qquad \text{for all } 0 < h < 1.
\end{equation}

2) In order that $I_\mu$ is bounded, it suffices that there is a constant $A > 0$ such that:
\begin{equation}\label{CS continue}
\qquad \quad K_{\mu, 2} (h) \leq \frac {1/h^2} {\Psi_2 [A \Psi_1^{-1} (1/h^2) ] } 
\, \raise 1,5pt \hbox{,} \qquad \text{for all } 0 < h < 1.
\end{equation}
\end{theoreme}

Note that condition \eqref{CN continue} reads as 
$\displaystyle \frac {\Psi_1^{-1} (1/h^2) } { \Psi_2^{-1} \big(1/ \rho_\mu (h) \big)}$ is bounded (by $1/ A$) and 
condition \eqref{CS continue} as 
$\displaystyle \frac {\Psi_1^{-1} (1/h^2) } { \Psi_2^{-1} \big(1/ h^2 K_{\mu, 2} (h) \big)}$ is bounded. 
\medskip

When $\Psi_1 = \Psi_2 = \Psi$ and the Orlicz function $\Psi$ satisfies the usual condition $\Delta_2$: 
$\Psi (2x) \leq C\, \Psi (x)$ for some constant $C > 1$ and $x$ large enough, it is clear that conditions 
\eqref{CN continue} and \eqref{CS continue} are equivalent. However, they are not equivalent in general; and even 
condition \eqref{CN continue} is not sufficient and condition \eqref{CS continue} is not necessary: 
the examples {\bf 1.b} and {\bf 2.} of \cite{CompOrli}, Chapter~4, \S \kern 1pt 3, given in 
the Hardy case, work also for the Bergman case. For the sake of completeness, we are going to sketch them.\par
\smallskip

\noindent{\bf Example 1.} {\it For every Orlicz function which does not satisfy the $\Delta_2$ condition, there 
exists a finite positive measure $\mu$ on $\D$ such that $I_\mu \colon {\mathfrak B}^\Psi \to L^\Psi (\mu)$ 
is continuous, though $\mu$ is not a $2$-Carleson measure, and {\it a fortiori} does not verify 
\eqref{CS continue}.}\par
\smallskip

\noindent{\bf Proof.} Since $\Psi$ does not satisfy $\Delta_2$, there exists an increasing sequence 
$(a_n)_{n \ge 1}$ such that $\Psi (2 a_n) /n$ is increasing and $\Psi (2 a_n) / \Psi (a_n) \ge n 2^n$. Define the 
discrete measure $\mu$:
\begin{displaymath}
\mu =\sum_{n=1}^\infty \bigg(\frac {n} {\Psi(2 a_n)} - \frac {n+1} {\Psi(2 a_{n+1})}\bigg)\delta_{x_n},
\end{displaymath}
where $x_n = 1 - 1/ \sqrt{\Psi (2 a_n)}$. As $\mu \big( [x_N, 1] \big) = N/ \Psi(2a_N)$, $\mu$ is not a 
$2$-Carleson measure. On the other hand, for every $f$ in the unit ball of ${\mathfrak B}^\Psi$, one has 
(\cite{CompOrli}, Lemma~5.2) $|f ( z) | \leq 8\, \Psi^{- 1} [ 1/ (1 - |z|)^2 ]$ and it is easy to check 
that, if $g (z) = \Psi^{-1} [ 1 / (1 - |z|)^2 ]$, then $\| g \|_{L^\Psi (\mu)} \leq 2$, so 
$\| f \|_{L^\Psi (\mu)} \leq 16$, proving that $I_\mu$ is bounded.
\qed
\par\smallskip

\noindent{\bf Example 2.} {\it Let $\Psi (x) = \e^x - 1$; there exists a finite positive measure $\mu$ on $\D$ 
such that \eqref{CN continue} holds but $I_\mu \colon {\mathfrak B}^\Psi \to L^\Psi (\mu)$ is not bounded.}
\par\smallskip

\noindent{\bf Proof.} Let $\nu$ be a probability measure on $\T$, supported by a compact set $L$ of Lebesgue 
measure zero, such that $\nu (I) \le |I|^{1/2}$, for each interval $I$. We can associate to $\nu$ the measure on 
$\overline{\D}$ defined by $\tilde\nu (E) =\nu (E\cap\T)$. By Rudin-Carleson's Theorem, for every integer 
$n$, there exists a function $g_n$ in the unit ball of the disk algebra such that $|g_n| = 1$ on $L$ and 
$\|g_n \|_{H^\Psi} \leq 4^{- n}$. As $L$ is compact, there exists some $r_n \in (1/2, 1)$ 
such that $|g_n (r_n z)| \ge 1/2$ for every $z\in L$. Now, define the measure $\mu$ by:
\begin{displaymath}
\mu (E) =\sum_{n=1}^\infty \frac {1} {2^n} \nu_n(E) \,,
\end{displaymath}
where:
\begin{displaymath}
\nu_n (E) =\nu \big(\{z\in\T \,;\  r_n z \in E\}\big).
\end{displaymath}
If $W$ is a Carleson window of size $h$ then, for each $n \ge 1$, we have:
\begin{displaymath}
\nu \big( \{z\in\T \,;\  r_n z\in W\} \big) \le \nu (W \cap \T) \le (2h)^{1/2}. 
\end{displaymath}
Hence, $\mu (W) \le (2h)^{1/2} \lesssim 1/ \Psi[ \frac{1}{4} \Psi^{-1} (1/h^2) ]$, and the condition 
\eqref{CN continue} is  fulfilled.\par

Nevertheless, the identity from  ${\mathfrak B}^\Psi$ to $L^1 (\mu)$ is not continuous since this would imply that 
the identity from  $H^\Psi$ to $L^1 (\mu)$ were continuous as well, which is not the case: 
$\| g_n \|_{L^1 (\mu)} \geq 1/ 2^{n + 1}$.
\qed

\par\bigskip

In order to prove Theorem~\ref{Theorem boundedness}, we shall need some results. They are analogous to 
Proposition~4.9, Theorem~4.13 and Lemma~4.14 of \cite{CompOrli}, but their proofs require different 
arguments\footnote{
By the way, we seize the opportunity to correct here the proof of Theorem~4.13 given in \cite{CompOrli}, where 
some argument had been put awkwardly. In that proof, we first had to set 
${\cal M} = \{ z \in \overline{\D}\,;\ |z| > 1 - h \quad \text{and} \quad |f (z) | > t\}$. Then, $M_f$ being 
the non-tangential maximal function of $f \in H^1$, the open set $\{M_f > t\}$ is the disjoint union of a 
countable family of open arcs $I_j \subseteq \T$, and we had to say that every $z$ such that $| f (z) | > t$ 
belongs to some window $W (I_j)$ (see \cite{Co-McC}, page 39). }.
\par\medskip

We first introduce the following maximal function:
\begin{equation}
\Lambda_f = \sum_{k = 0}^\infty \Big( \sup_{\Delta_k} |f| \Big)\,\ind_{\Delta_k}\,.
\end{equation}

One has:
\begin{lemme}\label{continuite fonction max}
For every Orlicz function $\Psi$, the map $f \in {\mathfrak B}^\Psi \mapsto \Lambda_f \in L^\Psi (\D)$ is 
bounded.
\end{lemme}

\noindent{\bf Proof.} Fix $f \in {\mathfrak B}^\Psi$. Set $c_k = \sup_{\Delta_k} |f|$ for every $k \geq 1$, and 
let $\alpha_k \in \Delta_k$ be such that $|f (\alpha_k) | \geq \frac{1}{2}\,\sup_{\Delta_k} |f| = c_k/2$. With 
$ C = \| f\|_{{\mathfrak B}^\Psi} > 0$, one has:
\begin{displaymath}
\int_\D \Psi (\Lambda_f/ 2C)\, d{\cal A} = \sum_{k \geq 0} \Psi (c_k /2 C) \, {\cal A} (\Delta_k) 
\leq \int_\D \Psi (|f| /C)\, d\mu\,,
\end{displaymath}
where $\mu = \sum_{k \geq 0} {\cal A} (\Delta_k) \,\delta_{\alpha_k}$.\par
But, for every Carleson window $W$, we can write:
\begin{displaymath}
\mu (W) = \sum_{\alpha_k \in W} {\cal A} (\Delta_k) 
\leq \sum_{\Delta_k \cap W \not= \emptyset} {\cal A} (\Delta_k) 
= {\cal A} \Big( \bigcup_{\Delta_k \cap W \not = \emptyset} \Delta_k \Big) \,,
\end{displaymath}
and, since $\bigcup_{\Delta_k \cap W \not= \emptyset} \Delta_k$ is contained in the window $\tilde W$ 
with the same center as $W$, but with size two times that of $W$, one has 
$\mu (W) \leq {\cal A} (\tilde W) = 4 \, {\cal A} (W)$. Hence $\mu$ is a Bergman-Carleson measure. By 
\cite{Hastings}, it follows that, for some constant $C_0 > 0$, which does not depend on $f$, one has, using the 
subharmonicity of $\Psi (|f|/C)$:
\begin{displaymath}
\int_\D \Psi (|f|/ C) \, d\mu \leq C_0 \int_\D \Psi (|f|/ C) \, d{\cal A} \leq C_0\,.
\end{displaymath}
We shall, as we may, assume that $C_0 \geq 1$. Now, by convexity of $\Psi$, we get:
\begin{displaymath}
\int_\D \Psi \Big( \frac {\Lambda_f} {2 C_0\, \| f \|_{{\mathfrak B}^\Psi} } \Big) \, d{\cal A} \leq 
\int_\D \frac{1}{C_0} \, \Psi \Big( \frac {\ \Lambda_f } {2 \, \|f \|_{{\mathfrak B}^\Psi} } \Big) \, d{\cal A} 
\leq 1 \,,
\end{displaymath}
meaning that $\| \Lambda_f \|_{L^\Psi (\D)} \leq 2 C_0\, \| f \|_{{\mathfrak B}^\Psi}$.
\qed
\medskip

\begin{lemme}\label{Carleson}
For every $f \in {\mathfrak B}^1$ and every finite positive Borel measure $\mu$ on $\D$, one has, for 
$0 < h < 1/2$ and $t > 0$:
\begin{displaymath}
\mu (\{z \in \D\,; \ |z| > 1 - h \quad \text{and} \quad |f (z) | > t\} ) 
\leq 4\, K_{\mu, 2} (2h) \,{\cal A} (\{\Lambda_f > t \}) \,.
\end{displaymath}
\end{lemme}

\noindent{\bf Proof.} Remark that when $z \in \Delta_k$ and $|z| > 1 - h$,  we must have 
$1 - 2^{- n - 1} > |z| > 1 - h$, hence $h > 2^{- n - 1}$; since $k = 2^n + j  - 1 \geq 2^{n - 1}$, we must have 
$k \geq 1/4h$. Let $I = \{ k \geq 1 \,; \ \sup_{\Delta_k} |f| > t \}$ and 
$I_h = \{ k \geq 1/4h\,; \ \sup_{\Delta_k} |f| > t \}$. If $W_k$ is the smallest Carleson window containing 
$\Delta_k$, we have:
\begin{align*}
\mu (\{z \in \D\,; \ |z| > 1 - h & \quad \text{and} \quad |f (z) | > t\} ) \\
& \leq \sum_{k \in I_h} \mu (\Delta_k)  \leq \sum_{k \in I_h, \, k \geq 1/4h} \mu (W_k)  \\
& \lesssim \sum_{k \in I_h} K_{\mu, 2} (2h) \, {\cal A} (W_k) 
\leq 4 \sum_{k \in I_h} K_{\mu, 2} (2h) \, {\cal A} (\Delta_k) \\
& \leq 4 K_{\mu, 2} (2h) \sum_{k \in I} {\cal A} (\Delta_k) 
= 4 K_{\mu, 2} (2h) \, {\cal A} (\{ \Lambda_f > t\}) \,. 
\end{align*}
and Lemma~\ref{Carleson} is proved. 
\qed

\begin{lemme}\label{L.U.}
Let $\mu$ be a finite Borel measure on $\D$ and $\Psi_1$ and $\Psi_2$ two Orlicz functions.\par
We suppose that, for some positive constant $A$, there is $0 < h_A \leq 1/2$ such that
\begin{displaymath}
\qquad \qquad K_{\mu, 2} (h) \leq \frac {1/h^2} {\Psi_2 [ A \Psi_1^{-1} (1/h^2) ] } \, \raise 1,5 pt \hbox{,} 
\quad \text{for} \quad 0 < h < h_A.
\end{displaymath}
Then, for every $f \in {\mathfrak B}^{\Psi_1}$ such that $\| f \|_{{\mathfrak B}^{\Psi_1}} \leq 1$ and every 
Borel subset $E$ of $\D$, one has, with $x_A = (A/8) \Psi_1^{-1} (4/h_A^2)$:
\begin{displaymath}
\int_E \Psi_2 (A\,|f| / 64) \, d\mu 
\leq \mu (E) \Psi_2 (x_A) + \frac {1}{8} \int_\D \Psi_1 (\Lambda_f)\,d{\cal A} \,.
\end{displaymath}
\end{lemme}

\noindent{\bf Proof.} For every $s > 0$, the inequality $| f (z)| > s$ implies that the norm of the evaluation 
$\delta_z$ at $z$ is greater than $s$. But this norm is (\cite{CompOrli}, Lemma~5.2) 
$\Psi_1 \big(1/ (1 - |z|)^2 \big)$, up to constants; more precisely:
\begin{displaymath}
\| \delta_z\| \leq 8 \, \Psi_1^{-1} \Big( \frac {1} {(1 - |z|)^2} \Big)\, \cdot
\end{displaymath}
Hence, we have:
\begin{displaymath}
s < 8 \, \Psi_1^{-1} \Big( \frac {1} {(1 - |z|)^2} \Big)\, \raise 1,5pt \hbox{,}
\end{displaymath}
so:
\begin{displaymath}
|z| > 1 - \frac {1} {\sqrt {\Psi_1 (s / 8) }}\, \cdot
\end{displaymath}
Lemma~\ref{Carleson} gives, when $\Psi_1 (s/8) \geq 2$:
\begin{align*}
\mu ( \{ |f (z ) | > s \}) 
& = \mu (\{ |z| > 1 - \frac {1} {\sqrt {\Psi_1 (s / 8) }} \quad \text{and} \quad|f (z ) | > s \} \\
& \leq 4\,K_{\mu, 2} \Big( \frac {2} {\sqrt {\Psi_1 (s / 8) }} \Big) \,{\cal A} (\{\Lambda_f > s \}) \,.
\end{align*}
But, by our assumption, if $\Psi_1 (s/8) \geq 4/h_A^2$,
\begin{align*}
K_{\mu, 2} \Big( \frac {2} {\sqrt {\Psi_1 (s / 8) }} \Big) 
& \leq \frac { \Psi_1 (s/8) /4} {\Psi_2 [ A \Psi_1^{-1} (\big( \Psi_1 (s/8) /4 \big) ] } \\
& \leq \frac { \Psi_1 (s/8) /4} {\Psi_2 [ (A/4) \Psi_1^{-1} (\big( \Psi_1 (s/8) \big) ] } 
\quad \text{(by convexity)} \\
& = \frac {1}{4}\,\frac { \Psi_1 (s/8) } {\Psi_2 (A s/32)} \, \raise 1,5 pt \hbox{;}
\end{align*}
hence:
\begin{displaymath}
\mu ( \{ |f (z ) | > s \}) \leq \frac { \Psi_1 (s/8) } {\Psi_2 (A s/32)} \,{\cal A} (\{ \Lambda_f > s \})\,.
\end{displaymath}

We get therefore:
\begin{align*}
\int_E \Psi_2 (A |f| / 64)\, d\mu 
& = \int_0^{+\infty} \Psi'_2 (t)\, \mu (\{ |f| > 64 t/A \} \cap E)\,dt \\
& \leq \int_0^{x_A} \Psi'_2 (t) \mu (E)\,dt \\
& \hskip 25 pt 
+ \int_{x_A}^{+\infty} \Psi'_2 (t) \, \frac {\Psi_1 (8 t/A)} {\Psi_2 (2t) } 
\, {\cal A} (\{ \Lambda_f > 64 t/ A\})\, dt \\
& \leq \Psi_2 (x_A) \mu (E) \\
& \hskip 25 pt 
+ \int_{x_A}^{+\infty} \frac {\Psi'_2 (t)} {\Psi_2 (2t)} \, \Psi_1 (8t/A) 
\, {\cal A} (\{ \Lambda_f > 64 t/ A\})\, dt\,.
\end{align*}

But, as $\Psi_1$ and $\Psi_2$ are Orlicz functions, one has $t \Psi'_2 (t) \leq \Psi_2 (2t)$ and 
$\Psi_1 (8 t/A) \leq (8t/A) \Psi'_1 (8t/A)$; hence:
\begin{align*}
\int_{x_A}^{+\infty} \frac {\Psi'_2 (t)} {\Psi_2 (2t)} \, \Psi_1 (8t/A)\, 
& {\cal A} (\{ \Lambda_f > 64 t/ A\})\, dt \\
& \leq \int_0^{+\infty} \frac {\Psi_1 (8t/A)} {t}\, {\cal A} (\{ \Lambda_f > 64 t/ A\})\, dt \\
& \leq \frac {8}{A} \int_0^{+\infty} \Psi'_1 (8t/A)\, {\cal A} (\{ \Lambda_f > 64 t/ A\})\, dt \\
& = \int_0^{+\infty} \Psi'_1 (x) \, {\cal A} (\{ \Lambda_f > 8x \})\, dx \\
& = \int_\D \Psi_1 (\Lambda_f /8) \, d{\cal A} \leq \frac {1}{8} \int_\D \Psi_1 (\Lambda_f)\, d{\cal A}\,,
\end{align*}
and the proof of Lemma~\ref{L.U.} is finished. 
\qed

\bigskip

\noindent{\bf Proof of Theorem~\ref{Theorem boundedness}.} 1) Consider, for every $a \in \D$, the Berezin 
kernel:
\begin{equation}\label{Ha}
H_a (z) = \frac {(1 - |a|^2)^2} {|1 - \overline{a} z|^4} \,\cdot 
\end{equation}
One has $\| H_a \|_{{\mathfrak B}^1} = 1$ and 
\begin{displaymath}
\| H_a \|_\infty = \frac {(1 - |a|^2)^2} {(1 - |a|)^4} = \frac {(1 + |a|)^2} {(1 - |a|)^2 } 
\leq \frac{4}{(1 - |a|)^2} \,;
\end{displaymath}
hence (\cite{CompOrli}, Lemma~3.9):
\begin{equation}\label{norme Psi de Ha}
\qquad \qquad \| H_a \|_{{\mathfrak B}^{\Psi_1}} \leq \frac {4/h^2} {\Psi_1^{-1} (4/h^2) } 
\, \raise 1,5 pt \hbox{,} \qquad \quad h = 1 - |a|\,.
\end{equation}
It follows that the function $f_a = \frac {1}{4}\, h^2 \Psi_1^{-1} (4/h^2) \, H_a$ is in the unit ball of 
${\mathfrak B}^{\Psi_1}$.\par
\smallskip

Now, let $\xi \in \T$ and $0 < h < 1$. When $z \in W (\xi, h)$, one has easily (see \cite{CompOrli}, proof of 
Theorem~4.10, with a slightly different definition of $W (\xi, h)$): $| 1 - \overline{a} z| \leq 5h$, where 
$a = (1 - h) \xi$. It follows that then $| f_a (z) | \geq (1/2500) \Psi_1^{-1} (4/h^2)$. Hence:
\begin{displaymath}
1 \geq \int_\D \Psi_2 \Big( \frac {|f_a|} {\| I_\mu \|} \Big)\,d\mu 
\geq \Psi_2 \Big( \frac {1}{2500\, \| I_\mu \|} \Psi_1^{-1} (4/h^2) \Big)\, \mu \big( W (\xi, h) \big) \,,
\end{displaymath}
which is \eqref{CN continue}.\par
\medskip

2) By Lemma~\ref{continuite fonction max}, there is a constant $C > 0$, that we may, and shall do, assume 
$\geq 1$, such that $\| \Lambda_f \|_{L^{\Psi_1} (\mu)} \leq C\, \| f \|_{{\mathfrak B}^{\Psi_1}}$ for every 
$f \in {\mathfrak B}^{\Psi_1}$. Let $g$ be in the unit ball of ${\mathfrak B}^{\Psi_1}$, and apply 
Lemma~\ref{L.U.} to $f = g/C$ (whose norm is $\leq 1$ yet), $E = \D$, with $h_A = 1/2$; we get, with 
$\tilde C = \max (1, \mu (\D) \Psi_2 (x_A) + \frac{1}{8} )$:
\begin{align*}
\int_\D \Psi_2 \Big( \frac {A} {64 C \tilde C} \, |g| \Big)\, d\mu 
& \leq \frac {1}{\tilde C} \int_\D \Psi_2 \Big( \frac {A} {64 C} \, |g| \Big)\, d\mu \\
& \leq \frac {1}{\tilde C} \, \bigg[ \mu (\D) \Psi_2 (x_A) 
+ \frac{1}{8} \int_\D \Psi_1 (\Lambda_f/C) \,d{\cal A} \bigg] \\
& \leq \frac {1}{\tilde C} \, \Big[ \mu (\D) \Psi_2 (x_A) + \frac{1}{8} \Big] \leq 1 \,,
\end{align*}
which means that $\| g\|_{L^{\Psi_2} (\mu)} \leq 64 C \tilde C /A$.
\qed


\subsection{Compactness}\label{compactness}

\begin{theoreme}\label{Theorem compactness}
Let $\mu$ be a finite positive Borel measure on $\D$, $\Psi_1$ and $\Psi_2$ two Orlicz functions, and let 
$I_\mu \colon {\mathfrak B}^{\Psi_1} \to L^{\Psi_2} (\mu)$ be the canonical map defined by $I_\mu (f) = f$. 
One has:
\par
1) If $I_\mu$ is compact, then: 
\begin{equation}\label{CN compact}
\lim_{h \to 0} \frac{\Psi_1^{-1} (1/ h^2)}{\Psi_2^{-1}\big( 1/\rho_\mu (h) \big)} = 0.
\end{equation}

2) In order that $I_\mu$ is compact, it suffices that 
\begin{equation}\label{CS compact}
\lim_{h \to 0} \frac{\Psi_1^{-1} (1/ h^2)}{\Psi_2^{-1}\big( 1/h^2 K_{\mu,2} (h) \big)} = 0.
\end{equation}
\end{theoreme}
\medskip

As for the boundedness case, conditions \eqref{CN compact} and \eqref{CS compact} are equivalent if 
$\Psi_1 = \Psi_2 = \Psi$ is sufficiently regular, but not in general. We shall give examples after the proof of the 
theorem, at the end of the section.
\smallskip
 
To prove the first part of this theorem, we shall need the following lemma.

\begin{lemme}\label{critere compact}
$I_\mu \colon {\mathfrak B}^{\Psi_1} \to L^{\Psi_2} (\mu)$ is compact if and only if for every bounded 
sequence $(f_n)_n$ in ${\mathfrak B}^{\Psi_1}$ converging to $0$ uniformly on compact subsets of $\D$, the 
sequence $\big( I_\mu (f_n) \big)_n$ converges to $0$ in the norm of $L^{\Psi_2} (\mu)$.
\end{lemme}

\noindent{\bf Proof of Lemma~\ref{critere compact}.} Assume that $I_\mu$ is compact, and let $(f_n)_n$
be a bounded sequence in ${\mathfrak B}^{\Psi_1}$ which converges to $0$ uniformly on compact subsets of 
$\D$. Since $I_\mu$ is compact and $(f_n)_n$ is bounded, we have a subsequence $(g_n)_n$ such that 
$g_n = I_\mu (g_n)$ converges to some $g \in L^{\Psi_2} (\mu)$. Then some subsequence of $(g_n)_n$ 
converges $\mu$-{\it a.e.} to $g$. Since $(g_n)_n$ converges to $0$ uniformly on compact subsets of $\D$, 
we must have $g = 0$ $\mu$-{\it a.e.}, that is $g = 0$ as an element of $L^{\Psi_2} (\mu)$. Now, the 
compactness of $I_\mu$ implies that $\| f_n \|_{L^{\Psi_2} (\mu)}$ tends to $0$.
\par
Conversely, assume that the condition of the lemma is satisfied, and let $(f_n)_n$ be an arbitrary bounded 
sequence in ${\mathfrak B}^{\Psi_1}$. Since the evaluation map is continuous on ${\mathfrak B}^{\Psi_1}$ 
(\cite{CompOrli}, Lemma~5.2), $(f_n)_n$ is a normal family, and Montel's Theorem gives a subsequence 
$(g_n)_n$ which converges uniformly on compact subsets, to some holomorphic function $g$. By 
Fatou's lemma, $g$ belongs to ${\mathfrak B}^{\Psi_1}$. Now, $(g_n - g)_n$ is a bounded sequence of 
${\mathfrak B}^{\Psi_1}$ which converges to $0$ uniformly on compact subsets of $\D$. By hypothesis, 
$\| g_n - g\|_{L^{\Psi_2} (\mu)}$ tends to $0$, and it follows that $I_\mu$ is compact. \qed
\par\bigskip

\noindent{\bf Proof of Theorem~\ref{Theorem compactness}.} 1) Assume that the map 
$I_\mu \colon {\mathfrak B}^{\Psi_1} \to L^{\Psi_2} (\mu)$ is compact. Consider, for every $a \in \D$, the 
Berezin kernel \eqref{Ha}. It follows from \eqref{Ha} and \eqref{norme Psi de Ha} that 
$H_a/ \|H_a\|_{{\mathfrak B}^{\Psi_1}}$ converges uniformly to $0$ on compact subsets of $\D$ as  
$|a|$ goes to $1$; hence, by Lemma~\ref{critere compact}, the compactness of $I_\mu$ implies that 
$\| (H_a / \| H_a \|_{{\mathfrak B}^{\Psi_1}}) \|_{L^{\Psi_2} (\mu)}$ tends to $0$. That means that, for 
every $\eps > 0$, one has:
\begin{displaymath}
\int_\D \Psi_2 \Big(\frac {\!\! |H_a|} { \eps \| H_a \|_{{\mathfrak B}^{\Psi_1}} } \Big)\, d\mu \leq 1 \,,
\end{displaymath}
for $|a|$ close enough to $1$, depending on $\eps$.\par

Now, let $\xi \in \T$ and $0 < h < 1$ with $a = (1 - h) \xi$. As already said in the proof of 
Theorem~\ref{Theorem boundedness}, $z \in W (\xi, h)$ implies that $|H_a (z)| \geq 1 / 625 h^2$. Therefore:
\begin{displaymath}
\int_\D \Psi_2 \Big(\frac { |H_a|} { \eps\, \| H_a \|_{{\mathfrak B}^{\Psi_1}} } \Big)\, d\mu 
\geq \Psi_2 \Big( \frac {1/625 h^2} {\eps (4/ h^2) / \Psi_1^{-1} (4/h^2) } \Big) \, \mu \big( W (\xi, h) \big).
\end{displaymath}
We get, for $h > 0$ small enough:
\begin{displaymath}
\mu \big( W (\xi, h) \big) \leq \frac {1} {\Psi_2 \big( (1/ 2500\eps) \Psi_1^{-1} (4/h^2) \big) } \, \cdot
\end{displaymath}
Since $\xi \in \T$ is arbitrary, it follows that, for $h > 0$ small enough:
\begin{displaymath}
\rho_\mu (h) \leq \frac {1} {\Psi_2 \big( (1/ 2500\eps) \Psi_1^{-1} (4/h^2) \big) } \, \raise 1,5pt \hbox{,}
\end{displaymath}
which reads:
\begin{displaymath}
\frac {\Psi_1^{-1} (4/h^2) } { \Psi_2^{-1} \big( 1/ \rho_\mu (h) \big) } \leq 2500 \eps.
\end{displaymath}

Since $\Psi_1^{-1} (4/h^2) \geq \Psi_1^{-1} (1/h^2)$, we have obtained \eqref{CN compact}. \qed
\par\medskip

2) Assume now that \eqref{CS compact} is satisfied. By Lemma~\ref{critere compact}, we have to show that 
for every sequence $(f_n)_n$ in the unit ball of ${\mathfrak B}^{\Psi_1}$ which converges uniformly to $0$ on 
compacts subsets of $\D$, $\big(I_\mu (f_n)\big)_n$ converges to $0$ for the norm of $L^{\Psi_2} (\mu)$. So, 
let $(f_n)_n$ be such a sequence, and let $\eps >0$. \par
By Lemma~\ref{continuite fonction max}, there is a constant $C \geq 1$ such that 
$\| \Lambda_f \|_{L^{\Psi_1} (\mu)} \leq C\, \| f\|_{{\mathfrak B}^{\Psi_1}}$ for every 
$f \in {\mathfrak B}^{\Psi_1}$. Set $A = 64 C/ \eps$. By \eqref{CS compact}, there is an $h_A < 1$ such that,  
for $0 < h \leq h_A$, one has $\Psi^{-1} (1/h^2) \leq (1/A)\, \Psi_2^{-1} \big(1/ h^2 K_{\mu, 2} (h) \big)$, 
{\it i.e.} 
\begin{displaymath}
K_{\mu, 2} (h) \leq  \frac {1/h^2} {\Psi_2 [A \Psi_1^{-1} (1/h^2) ] } \,\cdot
\end{displaymath}
For $0 < r < 1$, we may therefore apply, for every $n \geq 1$, Lemma~\ref{L.U.} to $f_n/C$ (which is in the 
unit ball of ${\mathfrak B}^{\Psi_1}$), with $E = \D \setminus r \overline{\D}$; we get:
\begin{align*}
\int_{\D \setminus r \overline{\D}} \Psi_2 (|f_n|/ \eps)\, d\mu 
& = \int_{\D \setminus r \overline{\D}} \Psi_2 (A |f_n|/ 64 C)\, d\mu \\
& \leq \mu (\D \setminus r \overline{\D}) \Psi_2 (x_A) 
+ \frac{1}{8} \int_\D \Psi_1 (\Lambda_{f_n}/ C) \,d {\cal A} \\
& \leq \mu (\D \setminus r \overline{\D}) \Psi_2 (x_A) + \frac{1}{8} \,\cdot
\end{align*}
But this last quantity is $\leq 1/2$ for $r$ small enough. \par
Fix such an $r < 1$. Since $(f_n)_n$ converges uniformly to $0$ on compacts subsets of $\D$, one has 
$\int_{r \overline{\D}} \Psi_2 (|f_n|/\eps) \,d\mu \leq 1/2$ for $n$ large enough.\par
It follows that $\int_\D  \Psi_2 (|f_n|/ \eps)\, d\mu \leq 1$, and hence 
$\| f_n \|_{L^{\Psi_2 (\mu)}} \leq \eps$, for $n$ large enough.\par
That ends the proof of Theorem~\ref{Theorem compactness}.
\qed
\bigskip

As we said, in general, condition \eqref{CN compact} is not sufficient to ensure compactness, and condition 
\eqref{CS compact} is not necessary. They are equivalent (and hence necessary and sufficient for compactness) if 
$\Psi_1 = \Psi_2 = \Psi$ and $\Psi$ is a \emph{regular} Orlicz function. Here, \emph{regular} means that $\Psi$ 
satisfies the condition we called $\nabla_0$: for some $x_0 > 0$ and some $C \geq 1$, one has 
$\frac {\Psi (2x)} {\Psi (x)} \leq \frac {\Psi (2C y)} {\Psi (y)}$ for $x_0 \leq x \leq y$ (see \cite{CompOrli}, 
Theorem~4.11, whose proof works as well in the Bergman case). However, we gave in \cite{CompOrli}, examples 
showing, in the  Hardy case,  that this is not always the case (examples {\bf 3} and {\bf 4} in 
\cite{CompOrli}, Chapter~4, \S \kern 1pt 3). These examples work in the Bergman case and we are going to recall 
them sketchily.\par
\smallskip

\noindent{\bf Example 1.} {\it For every Orlicz function $\Psi$ not satisfying $\nabla_0$, there exists a measure 
$\mu$ such that $I_\mu \colon {\mathfrak B}^\Psi \to L^\Psi (\mu)$ is compact but for which \eqref{CS compact} 
is not satisfied.}
\par\smallskip

\noindent{\bf Proof.}  Since $\Psi \notin \nabla_0$, we can select two increasing sequences 
$(x_n)_{n \ge 1}$ and $(y_n)_{n \ge 1}$, with $1 \le x_n \le y_n \le x_{n+1}$ and $\Psi (x_n) > 1$, such that 
$\lim x_n =+\infty$ and
\begin{displaymath}
\frac {\Psi (2 x_n)} {\Psi (x_n)} \ge \frac {\Psi (2^n y_n)} {\Psi (y_n)} \,\cdot
\end{displaymath}
Define the discrete measure
\begin{displaymath}
\mu =\sum_{n=1}^\infty \frac {1} {\Psi (2^n y_n)} \,\delta_{r_n} \,,
\end{displaymath}
where $r_n = 1 - 1/\sqrt{\Psi (y_n)}$. The series converge since $\Psi (2^n y_n) \ge 2^n$. \par
The same proof as in \cite{CompOrli} shows that $I_\mu$ is compact, but, writing 
$h_n = 1 / \sqrt{\Psi(x_n)}$ and $t_n = 1/ \sqrt{\Psi(y_n)}$, we have:
\begin{displaymath}
K_{\mu, 2} (h_n) \ge \frac {\mu ([1 - t_n, 1])} {t_n^2} \ge \frac {\Psi (y_n)} {\Psi (2^n y_n)} 
\ge \frac {\Psi (x_n)} {\Psi (2x_n)} =\frac {1/h_n^2} {\Psi \big( 2\Psi^{-1} (1/h_n^2) \big)} \,\raise 1,5pt\hbox{,}
\end{displaymath}  
showing that \eqref{CS compact} is not satisfied.
\qed
\medskip

This actually shows that condition $\Psi \in \nabla_0$ is necessary and sufficient in order to ensure that the identity 
from $B^\Psi$ to $L^\Psi(\mu)$ is compact if and only if $\mu$ satisfies (2.7). \par
\medskip

\noindent{\bf Example 2.} {\it There exist an Orlicz function $\Psi$ and a measure $\mu$ on $\D$ such that 
\eqref{CN compact} holds, but for which $I_\mu \colon {\mathfrak B}^\Psi \to L^\Psi (\mu)$ is not compact.}
\par\smallskip

\noindent{\bf Proof.} We shall use the Orlicz function $\Psi$ introduced in \cite{Critere}. The key properties of 
this function are: \par
1) $\Psi (x) \ge x^3/3$ for every $x > 0$; \par
2) $\Psi (k!) \le (k!)^3$ for every integer $k \ge 1$; \par
3) $\Psi \big( 3(k!) \big) > k.(k!)^3$ for every integer $k \ge 1$.\par
\smallskip

\noindent Define $x_k = k!$, $y_k =(k+1)! / k^{1/3}$, $r_k = 1 - 1/ \sqrt{\Psi (y_k)}$ and 
$\rho_k = 1 - 1/ \sqrt{\Psi(x_k)}$. Of course, $x_2 < y_2 < x_3 < \cdots$. Let $\nu$ be the discrete measure 
defined by: 
\begin{displaymath}
\nu = \sum_{k = 2}^\infty \nu_k \,,
\end{displaymath} 
where:
\begin{displaymath} 
\nu_k = \frac {1} {\Psi \big( (k+1)! \big)} \sum_{a^{k^2} = 1} \delta_{r_ka} \,.
\end{displaymath} 

In order to show that \eqref{CN compact} is satisfied, it is clearly sufficient to prove that, when 
$1/ \sqrt{\Psi (y_k)} \le h < 1/ \sqrt{\Psi (y_{k-1})}$ (with $k \ge 3$), we have:
\begin{displaymath} 
\rho_\nu (h) \le \frac {1} {\Psi \big( \frac{1}{2} k^{1/3} \Psi^{-1} (1/h^2) \big)} \,\cdot
\end{displaymath} 
But, for such $h$, we have $\Psi^{-1} (1/h^2) \le y_k$ so
\begin{displaymath} 
\Psi \big( {\textstyle \frac{1}{2}} k^{1/3} \Psi^{-1}(1/h^2) \big) \le \frac{1}{2} \,\Psi \big( (k+1)! \big).
\end{displaymath} 
Hence, the conclusion follows from the fact that $\rho_\nu (h) \le 2 / \Psi \big( (k+1)! \big)$ 
(see \cite{CompOrli} for more details). So, condition \eqref{CN compact} is fulfilled.\par
\smallskip

We now introduce $\displaystyle f_k (z) = x_k \Big(\frac {1 - \rho_k} {1 - \rho_k z^{k^2}} \Big)^4$\,. By 
\eqref{norme Psi de Ha}, $\| f_k \|_{{\mathfrak B}^\Psi} \leq 1$. An easy computation gives 
$r_k^{k^2} \ge \rho_k$, for every $k \ge 2$; so, for every $a\in\T$ with $a^{k^2} = 1$, we have:
\begin{displaymath} 
f_k (a r_k) \ge x_k \Big( \frac {1 - \rho_k} {1 - \rho_k^2} \Big)^4 \ge\frac {1} {16} \, x_k.
\end{displaymath} 
Hence:
\begin{align*}
\int_{\D \setminus r_{k - 1} \D} \Psi (48 \, |f_k|) \,d\nu
& \ge \int_{\D \setminus r_{k - 1} \D} \Psi (48 \, |f_k|) \,d\nu_k 
\ge \frac {k^2} {\Psi \big( (k+1)! \big)} \Psi (3 x_k) \\
& > \frac {k^2} {\Psi \big( (k+1)! \big)} \big( k.(k!)^3 \big) \ge 1.
\end{align*}
Therefore, we conclude that 
$\sup_{\|f\|_{{\mathfrak B}^\Psi} \le 1} \|f \|_{L^\Psi (\D \setminus r_k \D, \mu)} \ge 1/48$, though 
$r_k \to 1$. Hence (see the above proof of Theorem~\ref{Theorem compactness}), $I_\mu$ is not compact.
\qed


\section{Compactness for composition operators}\label{compactness composition op}

\subsection{Carleson function}

We know (\cite{CompOrli}, Proposition~5.4), that every analytic self-map $\phi \colon \D \to \D$ induces 
a bounded composition operator $C_\phi \colon {\mathfrak B}^\Psi \to {\mathfrak B}^\Psi$. The main 
result of this section is that, for the pull-back measure ${\cal A}_\phi$ of ${\cal A}$ by $\phi$, the 
necessary and the sufficient conditions of Theorem~\ref{Theorem compactness} are equivalent. The same 
kind of result occurs for Hardy-Orlicz spaces (\cite{CompOrli}, Theorem~4.18 and Theorem~4.19), but 
the proofs must be different (because we use the analytic functions themselves, and not their boundary values).
\par\medskip

We have the following contractivity (or homogeneity) result, which can be viewed as a ``multi-scaled'' fact that 
${\cal A}_\phi$ is a $2$-Carleson measure.
\begin{theoreme}\label{homogeneite}
There exists a constant $C_0 > 0$ such that, for every analytic self-map 
$\phi \colon \D \to \D$, one has:
\begin{equation}\label{majoration homogeneite}
{\cal A} \big( \{\phi \in S (\xi, \eps h) \} \big) 
\leq C_0 \,\eps^2 \, {\cal A} \big( \{\phi \in S (\xi, h) \} \big)
\end{equation}
for every $\xi \in \T$, $0 < h < (1 - | \phi (0)|)$, and $0 < \eps \leq 1$.
\end{theoreme}

As a consequence, one has $\rho_{\phi, 2} (\eps h) \leq C \, \eps^2\, \rho_{\phi, 2} (h)$, for $h > 0$ 
small enough, and hence:
\begin{displaymath}
\frac {\rho_{\phi, 2} (h)} {h^2} \leq K_{\mu, 2} (h) 
= \sup_{0 < \eps \leq 1} \frac {\rho_{\phi, 2} (\eps h)} {\eps^2 h^2} 
\leq C_\alpha\, \frac {\rho_{\phi, 2} (h)} {h^2} \,\cdot
\end{displaymath}
Therefore:
\begin{theoreme}\label{CNS compactness composition operators}
For every analytic self-map $\phi \colon \D \to \D$ and every Orlicz function $\Psi$, the composition 
operator $C_\phi \colon {\mathfrak B}^\Psi \to {\mathfrak B}^\Psi$ is compact if and only if
\begin{equation}
\lim_{h \to 0} \frac{\Psi^{-1} (1/ h^2)}{\Psi^{-1}\big( 1/\rho_{\phi, 2} (h) \big)} = 0.
\end{equation}
\end{theoreme}
\bigskip

\noindent{\bf Proof of Theorem~\ref{homogeneite}.} It suffices, even if it means enlarging $C_0$, to show 
\eqref{majoration homogeneite} for $0 < h \leq h_0 = \alpha_0 (1 - |\phi (0)|)$ and $0 < \eps \leq \eps_0$ 
for some $0 < \alpha_0 < 1$ and $0 < \eps_0 < 1$. Indeed, if $0 < h \leq h_0$ and $\eps_0 \leq \eps \leq 1$, 
one has 
${\cal A} \big( S (\xi, \eps h)\big) \leq {\cal A} \big( S (\xi,  h)\big) 
\leq (1/\eps_0^2)\,\eps^2 {\cal A} \big( S (\xi,  h)\big)$. Now, for $h_0 \leq h < 1 -|\phi (0)|$, one has, on 
the one hand, for $0 < \eps \leq \alpha_0$,  
${\cal A} \big( S (\xi, \eps h)\big) \leq {\cal A} \big( S (\xi, (\eps/ \alpha_0) h_0)\big) 
\leq C_0 (\eps^2/ \alpha_0^2)\, {\cal A} \big( S (\xi, h_0) \big) 
\leq (C_0/\alpha_0^2)\, \eps^2  {\cal A} \big( S (\xi, h) \big)$; and, on the other hand, for 
$\alpha_0 \leq \eps \leq 1$, 
${\cal A} \big( S (\xi, \eps h)\big) \leq {\cal A} \big( S (\xi,  h)\big) 
\leq (1/ \alpha_0^2) \, \eps^2 {\cal A} \big( S (\xi,  h)\big)$. \par
\smallskip

Since, moreover, it suffices to make the proof for $\xi =1$, Theorem~\ref{homogeneite} will result from 
the following theorem.
\begin{theoreme}\label{theo cle}
There exist a constant $K > 0$, $\alpha_0 > 0$ and $\lambda_0 > 1$ such that every 
analytic function $f \colon \D \to \Pi^+$ with $| f (0)| \leq \alpha_0$ satisfies, for every 
$\lambda \geq \lambda_0$:
\begin{displaymath}
{\cal A} (\{ | f| > \lambda\}) \leq \frac {K} {\lambda^2}\, {\cal A} (\{ |f| > 1\}) \,,
\end{displaymath}
where $\Pi^+$ is the right-half plane $\Pi^+ = \{z\in \C\,; \ \Re z > 0\}$.
\end{theoreme}
Indeed, let $f = h/ ( 1 - \phi)$. Then $\Re f > 0$ and $| f (0)| \leq h/ (1 - |\phi (0)|) \leq \alpha_0$. We 
may apply Theorem~\ref{theo cle} and we get, for $0 < \eps \leq 1/ \lambda_0$:
\begin{align*}
{\cal A}_\phi \big( S (1, \eps h) \big) 
& = {\cal A}\, (\{ |\phi - 1| < \eps h | \}) 
= {\cal A}\, (\{ |f| > 1/\eps\} ) \\
& \leq K \eps^2 {\cal A}\, (\{ |f|  > 1\}) 
= K \eps^2 {\cal A}\, (\{ |1  -\phi| < h \}) 
= K \eps^2 {\cal A}_\phi \big( S (1, h) \big)\,,
\end{align*}
which proves Theorem~\ref{homogeneite}.
\qed

\subsubsection{Some lemmas}

\begin{lemme}\label{lemme 1}
There is some constant $C_1 > 0$ such that
\begin{equation}\label{ineg 1}
{\cal A} ( \{|f| > \lambda \} ) \leq \frac {C_1} {\lambda^2} \, | f (0)|^2 \,,
\end{equation}
for every analytic function $f \colon \D \to \Pi^+$ and for every $ \lambda >0$.\par
In particular, there is a constant $K_1 > 0$ such that $\| f \|_{L^1 (\D)} \leq K_1 |f (0)|$ for every such a function.
\end{lemme}

\noindent{\bf Proof.}  We may assume that $|f (0)| = 1$.\par
The second assertion follows from the first one:
\begin{displaymath}
\int_\D |f|\, d{\cal A} = \int_0^{+\infty} {\cal A} ( \{|f| > \lambda \} ) \,d\lambda 
\leq \int_0^1 d\lambda + \int_1^{+\infty} \frac {C_1} {\lambda^2}\,d\lambda = 1 + C_1 := K_1 \,;
\end{displaymath}
To prove the first assertion, remark first that the left-hand side of inequality~\eqref{ineg 1} is $\leq 1$, so 
\eqref{ineg 1} is obvious for $\lambda \leq 2$ (with $C_1 \geq 4$). Assume that $\lambda > 2$, and 
set $\phi (z) = \frac {f (z) - f (0)} {f (z) + \overline{f (0)} }$\,. Then $| f (z)|  > \lambda$ implies 
$| \phi (z)  - 1 | = 2\, | \Re f (0)| / |f (z) + \overline{f (0)}| \leq 2 / (\lambda - 1) \leq 4/ \lambda$. But 
$\phi$ maps $\D$ into itself and hence induces a bounded composition operator 
$C_\phi \colon {\mathfrak B}^2 \to {\mathfrak B}^2$. It follows from \cite{Hastings} that ${\cal A}_\phi$ is a 
Bergman-Carleson measure, and hence (see the proof of Theorem~\ref{Theorem boundedness}, 1), with 
$\Psi_1 (x) = \Psi_2 (x) = x^2$):
\begin{align*}
{\cal A} (\{ |f| > \lambda \}) 
& \leq {\cal A}_\phi \big( S (1, 4/ \lambda) \big) 
\leq {\cal A}_\phi \big( W (1, 20/ \lambda) \big) \\
& \leq C_0\, \| C_\phi \|^2 / (\lambda/20)^2 = C'_0\, \| C_\phi \|^2/ \lambda^2 \,, 
\end{align*}
for some constant $C_0 \leq 2500^2/4$.
But $\| C_\phi \| \leq \frac {1 + |\phi (0)| } {1 - |\phi (0)|} = 1$ (\cite{Zhu}, Theorem~11.6, page 308), and the 
result follows.
\qed
\par\bigskip

Let: 
\begin{equation}
G = \{ z \in \C\,; \ |\arg z | < \pi/4 \} \,.
\end{equation}
By applying Lemma~\ref{lemme 1} to $f^2$, we get:

\begin{lemme}\label{lemme 2}
There exists a constant $C_2 > 0$ such that for every analytic function $f \colon \D \to G$, one has, 
for every $\lambda > 0$
\begin{equation}
{\cal A} (\{ |f| > \lambda \}) \leq \frac {C_2} {\lambda^4} \, |f (0)|^4 \,.
\end{equation}
\par

In particular there is a constant $K_2 > 0$ such that $\| f \|_{L^2 (\D)} \leq K_2 |f (0)|$ for every such a function.
\end{lemme}

Now, we shall have to replace $\D$ by some conformal copies of $\D$. This will be possible by the following 
version of Lemma~\ref{lemme 2}.

\begin{lemme}\label{lemme 3}
Let $\Omega$ be a bounded Jordan domain bounded by a ${\cal C}^1$ Jordan curve $J$ and let 
$h \colon \D\to \Omega$ be a Riemann map which extends to a bi-Lipschitz homeomorphism $h$ of 
$\overline{\D}$ to $\overline{\Omega}$ such that: 
\begin{displaymath}
\qquad a \leq \vert h' (z)\vert \leq b\,, \qquad \forall z\in \D \,.
\end{displaymath}
Then, there exists a constant $C > 0$, depending only on $b/a$, such that, for any $\lambda > 0$ and 
any analytic function $f \colon \Omega \to G$, one has: 
\begin{equation}
\mathcal{A} (\{\vert f\vert>\lambda\}) 
\leq \frac {C} {\lambda^4} \,\mathcal{A}(\Omega) \, |f (c)|^4\,,
\end{equation}
where $c = h (0)$.
\par
In particular, there is a constant $K = K (\Omega, c) > 0$ such that 
$\| f  \|_2 \leq K \,{\cal A} (\Omega) \, |f (c)|$ for all such functions.\par

Moreover $C > 0$ can be taken so that, for every positive harmonic function $u \colon \Omega \to \R_+$, 
one has:
\begin{equation}\label{harmonicite deformee}
\frac {1} {C}\, u (c) \leq \frac {1} {{\cal A} (\Omega)} \int_\Omega u \,d{\cal A} \leq C\, u (c)\,.
\end{equation}
\end{lemme}
 
\noindent{\bf Proof.} By the change-of-variable formula,  we have, setting $F = f \circ h$ and $w = h (z)$:
\begin{align*}
\mathcal{A} (\{\vert f\vert >\lambda\}) 
& =\int _{\Omega} \ind_{\{\vert f \vert>\lambda\}} (w)\, d\mathcal{A} (w) 
=\int_{\D} \ind_{\{\vert F\vert>\lambda\}} (z) \vert h'(z)\vert^2 \, d\mathcal{A} (z) \\
& \leq b^2 \int_{\D} \ind_{\{\vert F \vert >\lambda\}} (z)\, d\mathcal{A} (z) 
= b^2 \mathcal{A} (\{\vert F\vert >\lambda\}) \leq b^2\, \frac {C_2} {\lambda^4} \, |f (c)|^4\,.
\end{align*}
We implicitly used the fact that $F$ maps $\D$ to $G$ and that $\vert F (0) \vert = \vert f (c) \vert$, 
so we were allowed to use the previous lemma. Moreover, we have
\begin{displaymath}
\mathcal{A} (\Omega) =\int_{\Omega} d\mathcal{A}(w) =\int_{\D}\vert h'(z) \vert^2 \,d\mathcal{A}(z) 
\geq a^2 \mathcal{A}(\D) =a^2 \,; 
\end{displaymath}
so that finally:
\begin{displaymath}
\mathcal{A} (\{\vert f \vert >\lambda\}) 
\leq C_2\, \frac {b^2} {a^2}\, \frac {1} {\lambda^4} \, \mathcal{A}(\Omega) \, |f (c)|^4 
\mathop{ = }\limits^{def} \frac {C} {\lambda^4}\, \mathcal{A}(\Omega) \, |f (c)|^4\,.
\end{displaymath}
This ends the proof of the first part of Lemma~\ref{lemme 3} for $\Omega$, with $C = C_2 b^2/a^2$. \par
\smallskip

Finally, one has:
\begin{align*}
u (c) 
& = (u \circ h) (0)  = \int_\D (u \circ h) (z)\, d{\cal A} (z) \\
& \geq \frac{1}{b^2} \int_\D u \big( h (z)\big) \, |h' (z)|^2\, d{\cal A} (z) 
= \frac{1}{b^2} \int_\Omega u (w)  \,d{\cal A} (w) \,,
\end{align*}
which gives the right-hand side of \eqref{harmonicite deformee}, since 
${\cal A} (\Omega) \geq a^2$. The left-hand side is proved in the same way, since ${\cal A} (\Omega) \leq b^2$.
\qed
\bigskip

The last lemma is of a different kind.

\begin{lemme}\label{lemme 4}
For every analytic function $f \colon D (z_0, r) \to G$, one has:
\begin{displaymath}
{\cal A} (\{  \Re f  > \Re f (z_0) /2 \}) 
\geq \frac {1} {8 K_2^2} \, {\cal A} \big( D (z_0, r) \big) \,,
\end{displaymath}
where $K_2$ is the constant given by Lemma~\ref{lemme 2}.
\end{lemme}

\noindent{\bf Proof.} Recall Paley-Zygmund's inequality (see \cite{LQ}, Proposition~III.3), usually stated 
in probabilistic language: for any positive random variable $X$ on some probability space, one has, for 
$0 < a < 1$: 
\begin{equation}\label{PZ}
\P \big( X > a\, \esp (X) \big) \geq (1 - a)^2\, \frac {[\esp (X)]^{2}} {\esp (X^2)} 
\, \raise 1,5pt \hbox{,}
\end{equation}
where $\esp$ stands for expectation.\par 
For our problem, we will  take as probability space the disk $D (z_0, r)$ equipped with the probability  measure 
$d\mathcal{A} / r^2$, and as a random variable $X = \Re f := u$. We get from \eqref{PZ} that, 
since $u (z_0) = \esp (X)$ by the mean value property of harmonic functions: 
\begin{displaymath}
\mathcal{A} (\{u > a\, u (z_0) \}) 
\geq (1 - a)^2  \, \frac {[u (z_0)]^2} {\esp (u^2)} \, r^2 \,.
\end{displaymath}
Now, observe that $\sqrt{2}\, \Re w \geq |w|$ when $w \in G$; hence the function 
$g (z) = f  (z - z_0)/ [\sqrt{2}\,u (z_0)]$, which maps $\D$ into $G$, satisfies $| g (0) | \leq 1$ and 
Lemma~\ref{lemme 2} gives $\esp (u^2) \leq \esp (|f|^2) = \| f \|_2^2 \leq 2 K_2^2 [u (z_0)]^2$. We get:
\begin{displaymath}
\mathcal{A} (\{u > a\, u (z_0) \}) \geq \frac {(1 - a)^2} {2K_2^2}\, r^2 \,,
\end{displaymath}
which gives the desired result by taking $ a = 1/2$.
\qed


\subsubsection{Proof of Theorem~\ref{theo cle}.} 

For technical reasons, we are going to work with functions with range in the set $G$. Proving 
Theorem~\ref{theo cle} amounts to prove:
\begin{proposition}\label{prop theo cle}
There exist constants $K' > 0$, $\alpha_1 > 0$ and $\lambda_1 > 1$ such that every analytic 
function $f \colon \D \to G$ with $| f ( 0)| \leq \alpha_1$ satisfies, for $\lambda \geq \lambda_1$:
\begin{displaymath}
{\cal A} ( \{ |f| > \lambda \}) \leq \frac{K'}{\lambda^4}\, {\cal A} ( \{ |f| > 1 \})\,.
\end{displaymath}
\end{proposition}

It will be useful to note that $\sqrt{2}\, \Re w \geq |w|$ when $w \in G$.\par
\medskip

\noindent{\bf Proof.} The idea of the proof is to split $F_\lambda = \{ |f| > \lambda \}$ (actually, 
$F_\lambda$ will be conditionned by the more regular set $\{M_d f > 1\}$) in parts on which 
we shall be able to apply Lemma~\ref{lemme 3}. In order to construct these parts, we shall use a 
Calder{\'o}n-Zygmund type decomposition adapted to the geometry of the unit disk. 
We are going to recall the principle of this decomposition for the convenience of the reader.\par 
\smallskip

Before that, we have to say that we shall begin to work not with the function $f$ as is the statement of 
Proposition~\ref{prop theo cle}, but with this function multiplied by a constant (as we shall specify at the end 
of the proof). Nevertheless, we shall denote this new function by $f$ yet.\par
Now, remark that for $f \colon \D \to G$ with $|f (0)| \leq \alpha < 1$, one has:
\begin{equation}\label{couronne}
|f (z) | > 1 \quad \Longrightarrow \quad |z| > \beta \,,
\end{equation}
for $\beta = (1 - \alpha^2)/ (1 + \alpha^2)$. In fact, setting $a = f (0)$, the function  
$g = \frac {f^2 - a^2} {f^2 + \overline{a}^2}$ maps $\D$ into itself and vanishes at $0$. By 
Schwarz's lemma, $|z| \leq \beta$ implies $| g (z)| \leq \beta$, and since 
$f^2 = \frac {\overline{a}^2 g + a^2} { 1 - g}$, we get:  
\begin{displaymath}
|f (z)|^2 \leq \frac {|a|^2 \beta + |a|^2} {1 - \beta} = \frac {|a|^2} {\alpha^2} \leq 1 \,.
\end{displaymath}

The goal of the remark \eqref{couronne} is that, in order to make the Calder{\'o}n-Zygmund decomposition, 
we have to avoid the center $0$ of $\D$, which is not covered by the sets $S \in {\mathcal S}$ defined below.\par
\smallskip

For convenience, we shall take $\beta = 1/2$ ({\it i.e.} $\alpha = 1/\sqrt 3$), so we shall be only concerned by the 
annulus $\Gamma = \{z \in \D\,;\ 1/2 \leq |z| < 1\}$.\par\smallskip

In order to perform the Calder{\'o}n-Zygmund decomposition, we have to make a splitting of the annulus 
$\Gamma$. For that, consider the rectangle 
\begin{displaymath}
R_0 = \{x +iy \,;\ \log (1/2) \leq x < 0 \text{ and } 0 \leq y < 2\pi \}
\end{displaymath}
and the family 
${\mathcal R} = \bigcup_{n\in\N} \mathcal{R}_n$ of all dyadic half-open rectangles, where ${\mathcal R}_n$ is 
the family of sets
\begin{displaymath}
R = \{ x + iy \in \C\,;\ x_j \leq x < x_{j + 1} \quad \text{and} \quad y_k \leq y < y_{k + 1} \}\,,
\end{displaymath}
where $x_j = (1 - j 2^{- n}) \log (1/2)$ and $y_k = 2 k\pi /2^n$, $0 \leq j, k \leq 2^n - 1$, $n \geq 0$.\par

We now maps the annulus $\Gamma$ onto the unit disk $\D$, using the exponential map 
$x + i y \mapsto z = \e^{x + i y}$. We get a family $\mathcal{S} =\bigcup_{n\in\N} \mathcal{S}_n$.\par
Note that the jacobian $\e^{2x}$ of the transformation is between $\beta^2 = 1/4$ and $1$; hence the area of a 
set in ${\mathcal S}_{n - 1}$ is less than $16$ times that of a set in ${\mathcal S}_n$.\par
\smallskip

Recall now the Calder{\'o}n-Zygmund decomposition. Consider the conditional expectation $\esp_n |f|$ of $|f|$, 
given the $\sigma$-algebra generated by $\mathcal{S}_n$, namely: 
\begin{displaymath}
(\esp_{n} |f|) (z) =\sum_{S \in\mathcal{S}_n} \Big( \frac {1} {\mathcal{A}(S)} \int_{S} |f| \,d{\cal A} \Big) 
\,\ind_{S}(z) \,.
\end{displaymath}
The dyadic maximal function $M_{d} f$ is defined by:
\begin{displaymath}
M_{d} f (z) =\sup_{n} (\esp_{n} |f|) (z) \,.
\end{displaymath} 
Note that:
\begin{displaymath}
M_d f (z) =\sup_{S_z} \frac {1} {{\cal A} (S_z)} \int_{S_z} |f| \,d{\cal A} \,,
\end{displaymath}
where the supremum is taken over all $S_z \in {\cal S}$ containing $z$.\par

By Lebesgue's differentiation Theorem (or by the martingale convergence Theorem), we know that 
$(\esp_{n} |f|) (z)$ converges to $|f (z)|$,  $\mathcal{A}$-almost everywhere, so $M_d f (z) \leq 1$ 
implies $|f (z)| \leq 1$ for almost all $z \in\D$, and we can write: 
\begin{equation}\label{inclusion}
F_1 =\{|f| > 1 \} \subseteq \{ M_{d} f > 1 \} \cup N  \,, 
\end{equation}
where $N$ is a negligible set. Hence $F_\lambda \subseteq \{ M_{d} f > 1 \} \cup N$ for 
$\lambda \geq 1$. Now, the set $\{ M_{d} f > 1 \}$ can be decomposed in a disjoint union 
$\{ M_{d} f > 1 \} = \bigsqcup_{n} E_n$, where 
\begin{displaymath}
E_n = \{z \in \D \,; \ (\esp_{n} |f|) (z) > 1 \text{ and } (\esp_{j} |f|) (z)\leq 1 \text{ if } j < n\}.
\end{displaymath}
Since $\esp_n |f|$ is constant on the sets $S \in {\cal S}_n$, each $E_n$ can be in its turn decomposed 
into a disjoint union $E_n = \bigsqcup_{k} S _{n, k}$, where $S_{n, k} \in \mathcal{S}_n$.\par
By definition, for $z \in E_n$, one has $(\esp_n |f|) (z) \geq 1$ and hence 
\begin{displaymath}
\qquad \frac {1} {{\cal A} (S_{n, k})} \int_{S_{n, k}} |f| \, d{\cal A} \geq 1 \qquad \text{for } z \in S_{n, k} \,.
\end{displaymath}
But, on the other hand, $(\esp_{n - 1} |f|) (z) \leq 1$, and we have, if $z \in S_{n, k}$:
\begin{align*}
(\esp_n |f|) (z) 
& = \frac {1} {{\cal A}(S_{n, k})} \int_{S_{n, k}} |f| \,d{\cal A} \\
& \leq \frac {1} {{\cal A}(S_{n, k})} \int_{S_{n - 1, j}} |f| \,d{\cal A} 
\leq 16\, \frac {1} {{\cal A}(S_{n - 1, j})} \int_{S_{n - 1, j}} |f| \,d{\cal A}  \leq 16\,,
\end{align*}
where $S_{n - 1, j}$ is the set of rank $(n - 1)$ containing $S_{n, k}$.\par
\smallskip

Finally, reindexing the sets $S_{k, n}$,  we can write $\{ M_{d} f > 1 \}$ as a disjoint union 
\begin{equation}\label{union disjointe}
\{ M_{d} f > 1 \} = \bigsqcup_{l \geq 1} S_l \,, 
\end{equation}
for which:
\begin{equation}\label{encadrement moyennes}
1 \leq \frac {1} {\mathcal{A} (S_l)} \int_{S_l} |f| \,d{\cal A} \leq 16 \,,
\end{equation}

Equations \eqref{inclusion}, \eqref{union disjointe} and \eqref{encadrement moyennes} define the 
Calder{\'o}n-Zygmund decomposition of the function $f$. \par
\medskip

In order to apply Lemma~\ref{lemme 3}, we have to control (from above and \emph{from below}) the 
values of $|f|$ at the ``center'' of the sets $S_l$. One might think of doing this by using 
\eqref{encadrement moyennes}, but it is not always possible  (for example, the function with positive real part 
$1/ z^2$ is not integrable on the square of vertices $0$, $1 +i$, $1- i$, $2$). Nevertheless, since our function is 
holomorphic on a domain bigger than $S_l$, we may enlarge $S_l$ in order to use 
\eqref{encadrement moyennes}.\par
\smallskip

\begin{figure}
\centering
\begin{minipage}[b]{0.45\linewidth}
\centering
\includegraphics[width=\linewidth]{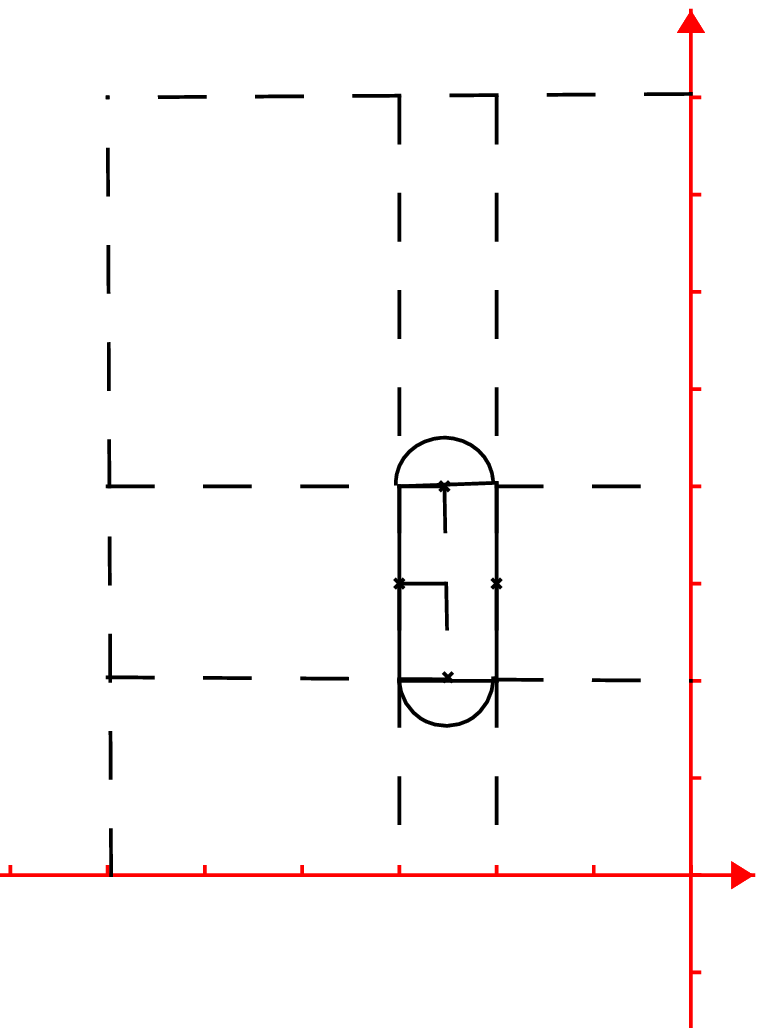} 
\caption{\it Round set $\hat R_l$} \label{rectangles}
\end{minipage}
\hspace{10pt}
\begin{minipage}[b]{0.45\linewidth}
\centering
\includegraphics[width=\linewidth]{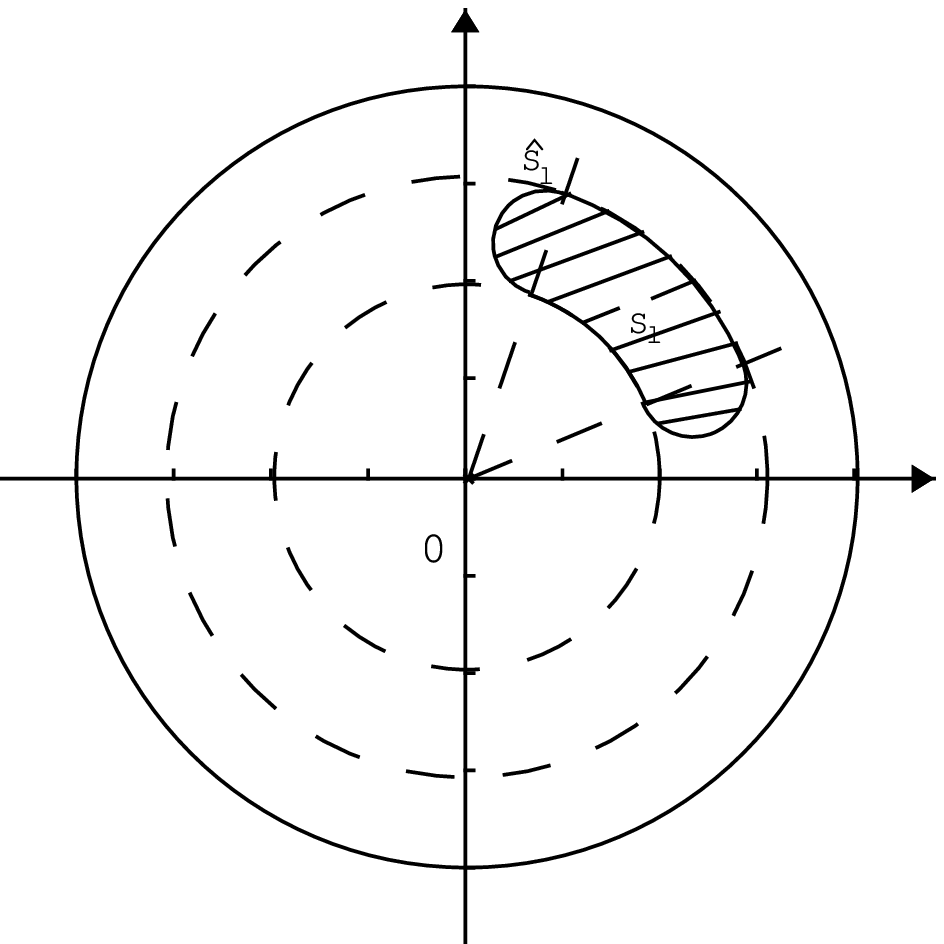} 
\vskip 20pt
\caption{\it Round set $\hat S_l$} \label{dyadic} 
\end{minipage}
\end{figure}

Let $R_l \in {\mathcal R}$ be the rectangle which is mapped onto $S_l$ by the exponential, and let us round 
$R_l$ in $\hat R_l$, by adding half-disks, as indicated in {\it Figure~\ref{rectangles}}. Let 
$\hat S_l = \exp (\hat R_l)$.

It is essential, when using Lemma~\ref{lemme 3} for $\Omega = \hat S_l$, that the constant $C$ given by this 
lemma does not depend on $l$. This will be checked in the following way. The sets $\hat R_l$ can be performed by 
making similarities from $\hat R_0$, where
\begin{displaymath}
R_0 = \{x +iy \,;\ \log (1/2) \leq x < 0 \text{ and } 0 \leq y < 2\pi \} \,.
\end{displaymath} 
\par

The boundary of $\hat R_0$ is ${\cal C}^1$ and it follows from \cite{Pom}, Theorem~3.5, for example, that 
$\hat R_0$ is  
conformally and bi-Lipschitz equivalent to $\D$; therefore, we are able to apply Lemma~\ref{lemme 3} (with 
$c = - \frac{1}{2} \log 2 + i \pi$). But if $T_l (\hat R_0) = \hat R_l$, with $T_l (z) = \alpha_l z + \beta_l$, 
$\alpha_l \neq 0$, we have, if $h_l = T_l \circ h$, 
$\vert {h_l} ' \vert = \vert \alpha_l \vert \vert h' \vert \leq \vert \alpha_l \vert \, b$, as well as 
$\vert {h_l} ' \vert \geq \vert \alpha_l \vert \, a$,  so that $a$ and $b$ are respectively changed into 
$a_l =\vert \alpha_l \vert a$ and $b_l = \vert \alpha_l \vert b$, and the quotient $b_l/a_l = b/a$ remains 
unchanged. Now, $\hat S_l = \exp (\hat R_l)$ and $1/2 \leq |\exp (x + iy)| = \e^x \leq 1$ for $x+iy \in \hat R_0$; 
hence the constant $C = C_2 b^2/a^2$ of Lemma~\ref{lemme 3} for $\Omega = \hat R_0$ becomes 
$C \leq 4 C_2 b^2 /a^2$ for $\Omega = \hat S_l$.\par
\smallskip

Therefore, one has, for every $l \geq 1$, if $\gamma_l$ is the center (defined in an obvious way) of $S_l$ (and 
hence of $\hat S_l$):
\begin{equation}
{\cal A} (\{ z \in \hat S_l \,;\ |f (z)| > \lambda \}) 
\leq \frac {C}{\lambda^4} \, |f (\gamma_l)|^4\, {\cal A} (\hat S_l) \,.
\end{equation}
\par\medskip

Now, let $D_l$ be the greatest open disk with center $\gamma_l$ and contained in $S_l$ (see 
{\it Figure~\ref{petit disque}}). We have, by the last part of Lemma~\ref{lemme 3}:

\begin{figure}
\centering
\includegraphics[width=0.3\linewidth]{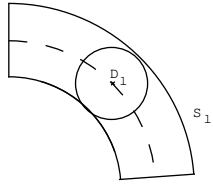} 
\caption{\it Disk $D_l$} \label{petit disque} 
\end{figure}

\begin{equation}
| f (\gamma_l) | 
\leq \frac {C} {{\cal A} (D_l)} \int_{D_l} |f | \, d{\cal A} 
\leq \frac {64 C} {{\cal A} (S_l)} \int_{S_l} |f| \, d{\cal A} \leq 2^{12}\,C \,.
\end{equation}
\par

Using the fact that ${\cal A} (\hat S_l) \leq 4\,{\cal A} (S_l)$ and that  
$F_\lambda \subseteq \{M_d f > 1 \} \cup N$ when $\lambda \geq 1$, we get:
\begin{align}
{\cal A} (F_\lambda) 
& = {\cal A} (F_\lambda \cap \{M_d f > 1\}) 
\leq \sum_{l \geq 1} {\cal A} (\{ z \in \hat S_l \,;\ |f (z)| > \lambda \})  \label{part 1}\\
&  \leq \frac {C}{\lambda^4} \, (2^{12} C)^4 \sum_{l \geq 1} {\cal A} (\hat S_l) 
\leq \frac {2^{52} C^5}{\lambda^4} \sum_{l \geq 1} {\cal A} (S_l) \notag \\
& = \frac {2^{52} C^5}{\lambda^4}\, {\cal A} ( \{ M_d f > 1\})\,. \notag
\end{align}

It remains to control ${\cal A} (\{ M_d f > 1\})$ by ${\cal A} (\{ |f | > \delta \})$, for some numerical 
$\delta > 0$. \par
For that, we shall use Lemma~\ref{lemme 4}. By harmonicity and Lemma~\ref{lemme 3}, one has, 
with $u = \Re f$:
\begin{align*}
u (\gamma_l) 
& \geq  \frac {1/C} {{\cal A} (\hat S_l)} \int_{\hat S_l} u \,d{\cal A} 
\geq  \frac {1} {16 C} \, \frac {1} {{\cal A} (S_l)} \int_{S_l} u \,d{\cal A} \\
& \geq \frac {1} {16\sqrt{2}\, C} \, \frac {1} {{\cal A} (S_l)} \int_{S_l} |f| \,d{\cal A} 
\geq \frac {1} {16\sqrt{2}\, C}  \,\cdot
\end{align*}

We now apply Lemma~\ref{lemme 4} and we get:
\begin{align*}
{\cal A} (\{ |f| > 1/ 64 C\} \cap D_l) 
& \geq {\cal A} (\{ u > 1/32 \sqrt{2}\, C \} \cap D_l)  
\geq {\cal A} (\{ u > u (\gamma_l) / 2\} \cap D_l) \\
& \geq \frac {1} {8 K_2^2} \, {\cal A} (D_l) \,.
\end{align*}
We obtain hence:
\begin{equation}\label{part 2}
{\cal A} (\{ M_d f > 1\}) 
 = \sum_{l \geq 1} {\cal A} (S_l) 
\leq 16 \sum_{l \geq 1} {\cal A} (D_l) 
\leq 128 K_2^2 \, {\cal A} (\{ |f| > 1/ 64 C\}) \,.
\end{equation}
\smallskip

The proof of Theorem~\ref{theo cle} is now finished, because \eqref{part 1} and \eqref{part 2} give, for 
$\lambda \geq 1$:
\begin{equation}\label{presque fini}
{\cal A} (\{ |f| > \lambda \}) 
\leq \frac {2^{38}\, C^5\, K_2^2} {\lambda^4}\, {\cal A} (\{ |f| > 1/64 C\}) 
\end{equation}
and if $f$ is as in the statement of Proposition~\ref{prop theo cle}, we can apply \eqref{presque fini} to 
$f_1 = 64 C\, f$ and we get, for $\lambda \geq \lambda_1 = 64 C$:
\begin{displaymath}
{\cal A} (\{ |f| > \lambda \}) 
\leq \frac {2^{38} (64 C)^4\, C^5  K_2^2} {\lambda^4}\, {\cal A} (\{ |f| > 1 \})\,,
\end{displaymath}
when $|f (0)| \leq 1/ 64 C$.
\qed


\subsection{Nevanlinna counting function}\label{Nevanlinna}

The Nevanlinna counting function is defined,  for every analytic function $\phi \colon \D \to \D$, and for 
every $w \in \phi(\D) \setminus \{\phi (0) \}$, by:
\begin{equation}
N_\phi (w) = \sum_{\phi (z) = w} \log \frac{1}{|z|} \,,
\end{equation}
where each term $\log \frac{1}{|z|}$ is repeated according to the multiplicity of $z$, and by 
$N_\phi (w) = 0$ for the other $w \in \D$. \par
Recall (see \cite{LLQR-N}) that, if $m$ is the normalized Lebesgue measure on $\T$, then the Carleson 
function of $\phi$ is the Carleson function of the pull-back measure $m_\phi$ of $m$ by $\phi$. We proved 
in \cite{LLQR-N} (Theorem~3.1 and Theorem~3.7) that the behaviour of $N_\phi$ is equivalent to that of the 
Carleson function $\rho_\phi$ in the following way.
\begin{theoreme}[\cite{LLQR-N}]\label{Nevanlinna-Carleson}
There exist two universal constants $C, c > 1$, such that, for every analytic self-map 
$\phi \colon \D \to \D$, one has:
\begin{equation}\label{equivalence Nevanlinna-Carleson}
\sup_{w \in W (\xi, h) \cap \D} N_\phi (w) \leq C\, \, m_\phi [ W (\xi, c\,h)] \,,
\end{equation} 
and
\begin{equation}\label{equivalence integrale N-C}
m_\phi \big( W (\xi, h) \big) \leq 
C\, \frac{1}{{\cal A} \big( W(\xi, c h) \big)} \, \int_{W(\xi, c h)} N_\phi (z)\,d{\cal A} (z)
\end{equation}
for $0 < h < 1$ small enough.
\end{theoreme}

We are going to deduce from Theorem~\ref{Nevanlinna-Carleson} the same result in the $2$-dim\-ensional 
case. The Nevanlinna counting function of order $2$ is defined (see \cite{Shapiro}, \S~6.2), for 
$w \in \D \setminus \{\phi (0)\}$, by:
\begin{equation}
N_{\phi, 2} (w) = \sum_{\phi (z) = w} \big( \log (1/ |z|) \big)^2 \,,
\end{equation}
where each preimage $z$ of $w$ appears as often as its multiplicity. The partial Nevanlinna counting function 
is defined, for $0 < r \leq 1$ by:
\begin{displaymath}
N_\phi (r, w) = \sum_{\phi (z) = w, |z| < r} \log (r / |z|) 
\end{displaymath}
and we have (\cite{Shapiro}, Proposition~6.6, where a misprint occurs):
\begin{equation}\label{formule integrale pour Nevanlinna}
N_{\phi, 2} (w) = 2 \int_0^1 N_\phi (r, w) \,\frac {dr} {r} \,\cdot
\end{equation}

One has:
\begin{theoreme}\label{equiv N-C Bergman}
There exists a universal constant $C > 1$, such that, for every analytic self-map 
$\phi \colon \D \to \D$, one has:
\begin{equation}
(1/C) \, \rho_{\phi, 2} (h/C) \leq \sup_{|w| \geq 1 - h} N_{\phi, 2} (w) \leq C\,\rho_{\phi, 2} (C\,h),
\end{equation} 
for $0 < h < 1$ small enough.
\end{theoreme}

\noindent{\bf Proof.} Let $w_0 = \phi (0)$ and set:
\begin{displaymath}
u (z) = \frac {w_0 - \phi (z) } {1 - \overline{w}_0 \phi (z) }\, \cdot
\end{displaymath}
Since $u (0) = 0$, Schwarz's lemma gives $| u (z) | \leq |z|$. Hence there is no $z$ with $|z| < t$ such that 
$\phi (z) = w$ when $t \leq | w_0 - w |/ |1 - \overline{w}_0 w| := |u_0 (w) | $. It follows that 
\eqref{formule integrale pour Nevanlinna} actually writes:
\begin{equation}\label{formule integrale 2}
N_{\phi, 2} (w) = 2 \int_{|u_0 (w) |}^1 N_\phi (r, w) \,\frac {dr} {r} \,\cdot
\end{equation}
\par\smallskip

1) It follows from \eqref{formule integrale 2} that:
\begin{displaymath}
N_{\phi, 2} (w) \leq \frac {2}{|u_0 (w) |^2} \int_0^1 N_\phi (r, w) \,r\,dr \,.
\end{displaymath}
But (see \cite{LLQR-N}, Lemma~3.4) $1/ |u_0 (w) | = |w - w_0| / |1 - \overline{w} w_0| > 1/3$ when 
$1 - |w| < (1 - |w_0|)/4$; therefore, for $1 - |w| < (1 - |w_0|)/4$, we have:
\begin{displaymath}
N_{\phi, 2} (w) \leq 18 \int_0^1 N_\phi (r, w) \,r\,dr \,.
\end{displaymath}
\par

Now, $N_\phi (r, w) = N_{\phi_r} (w)$, where $\phi_r (z) = \phi (r z)$, and 
it follows from \eqref{equivalence Nevanlinna-Carleson} that, for $w \in W (\xi, h)$, with $\xi\in \T$ and 
$h > 0$ small enough, one has:
\begin{align*}
N_{\phi, 2} (w) 
& \leq 18 \int_0^1 C m_{\phi_r}  [ W (\xi, c\,h)]\,r\,dr \\
& = 18 C \int_0^1 m \big(\{ \e^{i\theta} \,; \ \phi (r \e^{i\theta} ) \in W (\xi, c\,h) \} \big)\, r\,dr \\
& = 9 C {\cal A} \big( \{ z\in \D\,; \ \phi (z) \in W (\xi, c\,h) \} \big) \,.
\end{align*}
\par

2) Conversely, it follows from \eqref{formule integrale pour Nevanlinna} that:
\begin{displaymath}
N_{\phi, 2} (w) \geq 2 \int_0^1 N_\phi (r, w) \,r\,dr \,;
\end{displaymath}
hence:
\begin{align*}
\frac {1} {{\cal A} [W (\xi, ch)]} & \int_{ W (\xi, ch)} N_{\phi, 2} (z)\, d{\cal A} (z) \\
& \geq 2\int_0^1 \bigg[ \frac {1} {{\cal A} [W (\xi, ch)]} \int_{ W (\xi, ch)} N_\phi (r, z)\, d{\cal A} (z) \bigg]
\,r \,dr \\
& \geq \frac {2}{C} \int_0^1 m_{\phi_r} [W (\xi, h)]\, r\,dr \\
& = \frac {2}{C} \int_0^1 m \big( \{\e^{i\theta}\,;\ \phi (r \e^{i\theta}) \in W (\xi, h) \} \big) \,r\, dr \\
& = \frac{1}{C}\,{\cal A} \big( \{z \in \D \,;\ \phi (z) \in W (\xi, h) \} \big) \,,
\end{align*}
and that finishes the proof of Theorem~\ref{equiv N-C Bergman}.
\qed
\medskip

\noindent{\bf Remark.} Actually, the proof shows that for some constant $C > 1$, one has:
\begin{equation}\label{equiv mesure fenetre}
(1/ C) \, {\cal A}_\phi [W (\xi, h/C)] \leq \sup_{w \in W (\xi, h)} N_{\phi, 2} (w) 
\leq C\, {\cal A}_\phi [W (\xi, Ch)] \,,
\end{equation}
for every $\xi \in \T$ and $0 < h <1$ small enough. Since the $\ell_2$-norm is less than the $\ell_1$-norm, one 
has $N_{\phi, 2} (w ) \leq [N_\phi (w)]^2$, it follows hence from Theorem~\ref{Nevanlinna-Carleson} and 
Theorem~\ref{equiv N-C Bergman} (actually \eqref{equiv mesure fenetre}) that, for some constant 
$C > 1$, one has:
\begin{equation}\label{comparaison Carleson un et deux}
{\cal A} (\{ z\in \D\,;\ \phi (z) \in W(\xi, h) \}) 
\leq C\, \big[m (\{u\in \T\,;\ \phi^\ast (u) \in W (\xi, Ch)\}) \big]^2 \,,
\end{equation}
for every $\xi \in \T$ and every $0 < h < 1$ small enough, a fact which does not seem easy to proved in 
a straightforward way. In particular, for some constant $C >0$, one has, for $h > 0$ small enough:
\begin{displaymath}
\rho_{\phi, 2} (h) \leq C\, [\rho_\phi (Ch)]^2\,.
\end{displaymath}
\par

\begin{corollaire}\label{coroll Nevanlinna}
For every analytic self-map $\phi \colon \D \to \D$ and for every Orlicz function $\Psi$, the composition 
operator $C_\phi \colon {\mathfrak B}^\Psi \to {\mathfrak B}^\Psi$ is compact if and only if:
\begin{equation}\label{petit o Nevanlinna}
\lim_{h \to 0} \frac {\Psi^{-1} (1/h^2)} {\Psi^{-1} \big(1/ \nu_{\phi, 2} (h) \big) } = 0\,,
\end{equation}
where $\nu_{\phi, 2} (h) = \sup_{|w| \geq 1 - h} N_{\phi, 2} (w) $.
\end{corollaire}

\noindent{\bf Proof.} If $C_\phi$ is compact, Theorem~\ref{CNS compactness composition operators} gives, 
for every $A > 0$, an $h_A > 0$ such that, for $0 < h \leq h_A$:
\begin{displaymath}
\Psi^{-1} (1/h^2) \leq \frac {1} {A}\, \Psi^{-1} \big( 1/ \rho_{\phi, 2} (h) \big)\,.
\end{displaymath}
Then Theorem~\ref{equiv N-C Bergman} gives:
\begin{displaymath}
\nu_{\phi, 2} (h) \leq C\, \rho_{\phi, 2} (C h) \leq C / \Psi [ A \Psi^{-1} (1/ C^2 h^2) ] \,,
\end{displaymath}
{\it i.e.} $A\, \Psi^{-1} (1/ C^2 h^2) \leq \Psi^{-1} [C/ \nu_{\phi, 2} (h)]$; and then, by concavity, since $C > 1$:
\begin{displaymath}
\frac {1}{C^2}\, \Psi^{-1} (1/h^2) \leq \Psi^{-1} (1/C^2 h^2) 
\leq \frac {1} {A} \, \Psi^{-1} \big(C/ \nu_{\phi, 2} (h) \big) 
\leq \frac {C} {A} \, \Psi^{-1} \big(1/ \nu_{\phi, 2} (h) \big) \,,
\end{displaymath}
which implies \eqref{petit o Nevanlinna}.\par
The converse follows the same lines.
\qed

\begin{corollaire}\label{coro compacite}
The compactness of the composition operator $C_\phi \colon {\mathfrak B}^\Psi \to {\mathfrak B}^\Psi$ 
implies that:
\begin{displaymath}
\lim_{|z| \to 1} \frac {\Psi^{-1} \big( 1/ [1 - |\phi (z)|]^2 \big) } {\Psi^{-1} \big( 1/(1 - |z|)^2\big)} = 0 \,.
\end{displaymath}
\end{corollaire}

This corollary was proved in \cite{CompOrli}, Theorem~5.7, by a more direct method; and we also showed that 
the  condition is sufficient when $\Psi$ grows fast enough (namely, satisfies the condition $\Delta^2$). However, 
we do not know whether it is not sufficient for some symbol $\phi$ and some Orlicz function $\Psi$ (see the 
remark at the end of the paper). Nevertheless, it follows easily from Corollary~\ref{coroll Nevanlinna} that this 
condition is sufficient when $\phi$ is finitely-valent (see \cite{LLQR-N}, proof of Theorem~5.3).\par
\medskip

\noindent{\bf Proof.} Since 
$N_{\phi, 2} \big( \phi (z) \big) \geq \big( \log (1 /|z|) \big)^2 \geq ( 1 - |z|)^2$, it follows from 
Corollary~\ref{coroll Nevanlinna} that, for every $A > 0$, one has:
\begin{displaymath}
( 1 -|z|)^2 \leq \frac {1} {\Psi \big[ A \Psi^{-1} \big( 1/ ( 1 - |\phi (z)|)^2 \big) \big]} \,;
\end{displaymath}
that is:
\begin{displaymath}
\frac { \Psi^{-1} \big( 1/ ( 1 - |\phi (z)|)^2 \big) } {\Psi^{-1} \big(1/ (1 - |z|)^2\big)} \leq 1/ A\,,
\end{displaymath}
and that proves Corollary~\ref{coro compacite}.
\qed




\section{Comparison of the compactness of composition operators on Har\-dy-Orlicz spaces and on 
Bergman-Orlicz spaces} \label{Hardy-Bergman}

In the classical case ($\Psi (x) = x^p$, $1 \leq p < \infty$), it is known (\cite{McCluer-Shapiro}, 
Theorem~3.5, with Proposition~2.7) that the compactness of $C_\phi \colon H^p \to H^p$ implies the 
compactness of $C_\phi \colon {\mathfrak B}^p \to {\mathfrak B}^p$. On the other hand, we implicitly 
proved in \cite{CompOrli}, Theorem~5.7, that when $\Psi$ grows very fast (namely, satisfies the so-called 
$\Delta^2$ condition), then the compactness of $C_\phi \colon H^\Psi \to H^\Psi$ implies the 
compactness of $C_\phi \colon {\mathfrak B}^\Psi \to {\mathfrak B}^\Psi$. Let us write why: it is easy to show  
(see \cite{LLQR-N}, proof of Theorem~4.3) that the compactness of $C_\phi \colon H^\Psi \to H^\Psi$ 
implies that $\Psi^{-1} [1 / 1 - |\phi (z)|) ] / \Psi^{-1} [1/ (1 - |z|) ]$ tends to $0$ as $|z|$ goes to $1$, and 
we actually proved in \cite{CompOrli}, Theorem~5.7, that, when $\Psi \in \Delta^2$, this last condition implies 
the compactness of $C_\phi \colon {\mathfrak B}^\Psi \to {\mathfrak B}^\Psi$. The next proposition gives a 
condition on $\Psi$ which, though not very satisfactory, includes the cases $\Psi (x) = x^p$ and 
$\Psi \in \Delta^2$ and for which compactness on $H^\Psi$ implies compactness on ${\mathfrak B}^\Psi$. 

\begin{proposition}\label{Hardy donne Bergman}
Assume that the Orlicz function $\Psi$ satisfies the following condition: for every $A > 0$, there exist 
$x_A > 0$ and $B \geq A$ such that:
\begin{equation}\label{condition sur Psi}
\Psi [ A \Psi^{-1} (x^2) ] \leq \big( \Psi [B \Psi^{-1} (x)] \big)^2 
\end{equation}
for every $x \geq x_A$. Then every analytic map $\phi \colon \D \to \D$ defining a compact composition 
operator $C_\phi \colon H^\Psi \to H^\Psi$ also defines a compact composition operator 
$C_\phi \colon {\mathfrak B}^\Psi \to {\mathfrak B}^\Psi$.
\end{proposition}

It is clear that $\Psi$ satisfies \eqref{condition sur Psi} if $\Psi (x) = x^p$. Recall that condition $\Delta^2$ 
means that there exists some $\alpha > 1$ such that $[\Psi (t)]^2 \leq \Psi (\alpha t)$ for $t$ large enough. 
For such a function, one has $\Psi^{-1} (x^2) \leq \alpha \Psi^{-1} (x)$, and hence 
$\Psi [ A \Psi^{-1} (x^2) ] \leq \Psi [\alpha A \Psi^{-1} (x)]$, which is 
$\leq \big( \Psi [\alpha A \Psi^{-1} (x)] \big)^2$ for $x$ large enough, since it tends to infinity. 
Therefore, \eqref{condition sur Psi} holds if $\Psi$ has $\Delta^2$. Classical examples of Orlicz functions with 
$\Delta^2$ are $\Psi (x) = \e^{x^q} - 1$, where $q \geq 1$.\par
Another example, which does not have $\Delta^2$, but satisfies \eqref{condition sur Psi}, is 
$\Psi (x) = \exp \big( [\log (x +1)]^2 \big) - 1$. 
\medskip

\noindent{\bf Proof of Proposition~\ref{Hardy donne Bergman}.} One has   $N_{\phi, 2} (w) \leq [N_\phi (w)]^2$ 
since the $\ell_2$-norm is less than the $\ell_1$-norm. Let $A > 0$ be arbitrary. If $C_\phi$ is compact on 
$H^\Psi$, we know, by \cite{LLQR-N}, Theorem~4.2, that:
\begin{displaymath}
\sup_{|w| \geq 1 - h} N_\phi (w) = o\, \Bigg( \frac {1} {\Psi [B \Psi^{-1} (1/h)]} \bigg) \,,\quad 
\text{as } h \to 0.
\end{displaymath}
By \eqref{condition sur Psi}, we get:
\begin{displaymath}
\sup_{|w| \geq 1 - h} N_{\phi, 2} (w) = o\, \Bigg( \frac {1} {\Psi [A \Psi^{-1} (1/h^2)]} \bigg) \,\cdot 
\end{displaymath}
Corollary~\ref{coroll Nevanlinna} ensures that $C_\phi$ is compact on ${\mathfrak B}^\Psi$.
\qed

\medskip

However, we are going to see that the conclusion of Proposition~\ref{Hardy donne Bergman} does not hold for an 
arbitrary Orlicz function, by proving the following theorem.

\begin{theoreme}\label{theo Hardy-Bergman}
There exist an analytic self-map $\phi \colon \D \to \D$ and an Orlicz function $\Psi$ such that 
$C_\phi \colon H^\Psi \to H^\Psi$ is compact whereas 
$C_\phi \colon {\mathfrak B}^\Psi \to {\mathfrak B}^\Psi$ is not compact.
\end{theoreme}

In order to prove it, we shall show and use the following result.

\begin{theoreme}\label{symbol}
There exists an analytic self-map $\phi \colon \D \to \D$ such that, for some constants $c_2 \geq \pi$ and 
$\pi/4 \geq c_1 > 0$, one has, for some constant $C > 0$ and for $h > 0$ small enough:
\begin{align}
\rho_\phi (h) & \leq C\,\e^{- c_1/h}\,;\\
\rho_{\phi, 2} (h) & \geq (1/C)\,\e^{-c_2/h}\,.
\end{align}
\end{theoreme}

\noindent{\bf Proof of Theorem~\ref{theo Hardy-Bergman}.} 
We know (\cite{CompOrli}, Theorem~4.18) that $C_\phi \colon H^\Psi \to H^\Psi$ is compact if and only if 
\begin{equation}
\lim_{h \to 0} \frac{\Psi^{-1} (1/h)}{\Psi^{-1} \big( 1/\rho_\phi (h) \big) } = 0
\end{equation}
and (by Theorem~\ref{CNS compactness composition operators}) that the compactness of 
$C_\phi \colon {\mathfrak B}^\Psi \to {\mathfrak B}^\Psi$ is equivalent to
\begin{equation}
\lim_{h \to 0} \frac{\Psi^{-1} (1/h^2)}{\Psi^{-1} \big( 1/\rho_{\phi, 2} (h) \big) } = 0\,.
\end{equation}

Hence, it suffices to construct an Orlicz function $\Psi$ such that:

\begin{equation}\label{condition 1 pour Psi}
\lim_{x \to \infty} \frac{\Psi^{-1} (x) } {\Psi^{-1} (\e^{\,c_1 x}) } = 0 \,; 
\end{equation}
and 
\begin{equation}\label{condition 2 pour Psi}
\limsup_{x \to \infty} \frac{\Psi^{-1} (x^2) } {\Psi^{-1} (\e^{\,c_2 x} ) } > 0\,.
\end{equation}

We shall actually construct an increasing concave function $f \colon [0, \infty) \to \R_+$ such that 
$f (0) = 0$ and $f(\infty) = \infty$  which satisfies \eqref{condition 1 pour Psi} and 
\eqref{condition 2 pour Psi} (with $f$ instead of $\Psi^{-1}$) and we shall then take $\Psi = f^{-1}$.\par

1) For that, we set $\alpha_0 = 0$ and we define an increasing sequence of positive numbers 
$\alpha_1 = 1, \alpha_2, \ldots$ by:
\begin{equation}
\qquad \quad \alpha_{n+1} = \e^{c_1 \alpha_n}\,,\qquad n\geq 0\,,
\end{equation}
and we take $f$ affine on each interval $[\alpha_n, \alpha_{n+1}]$, $n \geq 0$. More precisely, we set
\begin{equation}
\qquad \qquad f (t) = A_n t + B_n \,, \qquad 
\text{for \quad} \alpha_{n - 1} \leq t \leq \alpha_n\,, \quad n \geq 1,
\end{equation}
where $A_1 = 1$, $B_1 = 0$, and for $n \geq 1$:
\begin{equation}\label{raccord}
B_{n+1} - B_n = (A_n - A_{n+1} ) \,\alpha_n \,.
\end{equation}
and, for $n \geq 0$:
\begin{equation}\label{condition lim sup}
\frac{B_{n + 1}}{A_{n + 1}} = \frac{1}{2} \,( \e^{c_2 \sqrt{\alpha_n}} - 3 \alpha_n)\,.
\end{equation}
Condition \eqref{raccord} ensures that $f$ is continuous. It is clear that $f$ is increasing and that 
$f (\infty) = \infty$.\par

Now, since $c_2 > \sqrt 6$, the function $u$ defined by $u (x) = \e^{c_2 x} - 3 x^2$ is positive and 
increasing for $x > 0$; hence, if one sets
\begin{equation}
\beta_n = (\e^{c_2 \sqrt{\alpha_n}} - 3 \alpha_n)/2 \,, \quad n \geq 0,
\end{equation}
$(\beta_n)_n$ is an increasing sequence of positive numbers. But $\beta_n = B_{n + 1} / A_{n + 1}$; hence  
\eqref{raccord} gives, for $n \geq 1$:
\begin{displaymath}
A_n = \frac{\alpha_n + \beta_n}{\alpha_n + \beta_{n - 1}}\, A_{n + 1}\,\cdot
\end{displaymath}
Since $\beta_n > \beta_{n - 1}$, it follows that $A_n > A_{n +1}$, and so the function $f$ is concave.\par
\smallskip

2) For $n$ large enough, one has 
$\alpha_n < \e^{c_2 \sqrt{\alpha_n}} < \e^{c_1 \alpha_n} = \alpha_{n + 1}$; hence, for these $n$,
\begin{displaymath}
\frac { f ( \alpha_n) } {f (\e^{c_2 \sqrt{\alpha_n}} )}
= \frac { A_{n +1 } \alpha_n + B_{n + 1}} {A_{n + 1} \e^{c_2 \sqrt{\alpha_n}} + B_{n +1} } 
= \frac {\alpha_n + \beta_n} {\e^{c_2 \sqrt{\alpha_n}} + \beta_n} = \frac{1}{3} \,\raise 1,5pt \hbox{,}
\end{displaymath}
and it follows that
\begin{displaymath}
\limsup_{x \to \infty} \frac {f (x^2)} {f (\e^{c_2 x}) } \geq \frac{1}{3}\,\cdot
\end{displaymath}
Condition \eqref{condition 2 pour Psi} is satisfied.\par
\smallskip

3) It remains to check condition~\eqref{condition 1 pour Psi}.\par

For that, we shall fix a number $M > c_2/ c_1$ and take $n_0$ large enough to have 
$\alpha_{n - 1} < M \sqrt{\alpha_n} < \alpha_n$ for $n \geq n_0$.\par

Let $x_0$ be such that $x \geq x_0$ if and only if $\alpha_{n - 1}  \leq x < \alpha_n$ with $n \geq n_0$. 
Choose such an $x$. We have:
\begin{displaymath}
\alpha_n = \e^{c_1 \alpha_{n - 1} } \leq \e^{c_1 x} \leq \e^{c_1 \alpha_n} = \alpha_{n + 1}
\end{displaymath}

We shall separate two cases. For convenience, we set $\eps_n = 1/ \beta_n$.
\par\smallskip

a) \emph{Case 1}: $\alpha_{n - 1}  \leq x < M \sqrt {\alpha_n}$. Then:
\begin{align*}
\frac {f (x)} {f (\e^{c_1 x}) } 
& = \frac {A_n x + B_n} {A_{n + 1} \e^{c_1 x} + B_{n + 1} } 
\leq  \frac {A_n M \sqrt{\alpha_n} + B_n} {A_{n + 1} \alpha_n + B_{n + 1} } \\
& = \frac {A_n M \sqrt{\alpha_n} + B_n} {A_n  \alpha_n + B_n } 
\, \raise 1,5pt \hbox{,} \qquad \text{by \eqref{raccord}} \\
&  = \frac {\eps_{n - 1} M \sqrt{\alpha_n} + 1} {\eps_{n - 1} \alpha_n + 1 } \sim \frac{M}{\sqrt {\alpha_n}} 
\, \raise 1,5pt \hbox{,}
\end{align*}
since
\begin{displaymath}
\eps_{n - 1} \sqrt{\alpha_n} 
= 2\,\frac {\e^{c_1 \alpha_{n - 1} /2}} {\e^{c_2 \sqrt{\alpha_{n - 1}}} - 3 \alpha_{n - 1} } 
\sim 2 \exp \Big( \frac{c_1}{2} \alpha_{n - 1} - c_2 \sqrt{\alpha_{n - 1}} \Big) 
\mathop{\longrightarrow}\limits_{n \to \infty} +\infty.
\end{displaymath}
\par\smallskip

b) \emph{Case 2}: $M\sqrt{\alpha_n} \leq x < \alpha_n$. Then:
\begin{align*}
\frac {f (x)} {f (\e^{c_1 x}) } 
& \leq \frac {f (\alpha_n)} {f (\e^{c_1 M \sqrt{\alpha_n}} )} 
= \frac {A_{n + 1} \alpha_n + B_{n + 1}} {A_{n + 1} \e^{c_1 M \sqrt{\alpha_n}} + B_{n + 1}} 
= \frac {\eps_n \alpha_n + 1} {\eps_n \e^{c_1 M \sqrt{\alpha_n}} + 1} \\
& = \ \frac {\displaystyle \frac {2 \alpha_n} {\e^{c_2 \sqrt{\alpha_n}} - 3 \alpha_n} + 1} 
{\displaystyle \frac {2 \e^{c_1 M \sqrt{\alpha_n}}} {\e^{c_2 \sqrt{\alpha_n}} - 3\alpha_n} + 1} 
\sim \exp\big( (c_2 - M c_1) \sqrt{\alpha_n} \big) \mathop{\longrightarrow}_{n \to \infty} 0
\end{align*}
since $M c_1 > c_2$.\par

Putting the two cases together, we get that 
$\displaystyle \lim_{x \to +\infty} \frac {f (x)} {f (\e^{c_1 x}) } = 0$, so \eqref{condition 1 pour Psi} is 
satisfied, and Theorem~\ref{theo Hardy-Bergman} is fully proved. \qed 

\bigskip

\begin{figure}[ht]
\centering
\includegraphics[width=6cm]{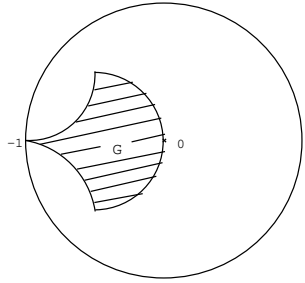}
\caption{\it Domain $G$} \label{domain}
\end{figure}

\noindent{\bf Proof of Theorem~\ref{symbol}.} The analytic map $\phi$ will be a conformal mapping from 
$\D$ to the domain $G$, edged by three circular arcs of radii $1/2$, and which is represented in 
{\it Figure~\ref{domain}}.\par
More precisely, let $G_0 = \D \cap \{ \Re z > 0\}$ and let $f \colon \D \to G_0$ be the conformal map 
such that 
\begin{displaymath}
f (- 1 ) = 0\,; \quad f (1) = 1 \,; \quad f (i) = i\,; \quad f (- i) = - i\,.
\end{displaymath}
We define successively $\phi_1 (z) = \log f (z)$, which maps $\D$ onto the half-strip  
$\{ \Re w < 0\, \ |\Im w| < \pi/2\}$, $\phi_2 (z) = - \frac{2}{\pi} \phi_1 (z) + 1$, 
$\phi_3 (z) = \frac{1}{\phi_2 (z)}$, and finally $\phi (z) = \phi_3 (z) - 1$.

1) When $W (\xi, h) \cap G \not= \emptyset$, we must have $W (\xi, h) \cap G \subseteq S (- 1, 2h)$, for 
$h$ small enough. Hence:
\begin{align*}
\rho_\phi (h) 
& \leq m (\{z \in \T\,;\ |\phi_3 (z) | < 2h\}) 
= m (\{z \in \T\,;\ |\phi_2 (z) | > 1/ 2h\}) \\
& \leq m (\{z \in \T\,;\ \Re \phi_2 (z)  > 1/ 2h - 1\})\,, \qquad \text{since } |\Im \phi_2 (z) | < 1\,, \\
& = m (\{z \in \T\,;\ \Re \phi_1 (z) < \pi - \pi/ (4h) \} ) \\
& = m (\{z \in \T\,; \ | f(z)| < \e^\pi \e^{-\pi/4h} \}) \\
& \leq K\, |\{t \in [-1, 1]\,;\ | it | < \e^\pi \e^{-\pi/4h} \}| = 2K \,\e^\pi \e^{-\pi/4h} \,,
\end{align*}
for small $h > 0$.\par
\smallskip

2) On the other hand, $S (-1, h) \cap G \subseteq W (-1, h)$, so:
\begin{align*}
\rho_{\phi, 2}  (h) 
& \geq {\cal A} (\{z \in \D\,;\ \phi (z) \in W (-1, h)\}) \\
& \geq {\cal A} (\{z \in \D\,;\ |\phi_3 (z) | < h\} ) 
= {\cal A} (\{z \in \D\,;\ |\phi_2 (z) | > 1/ h\} ) \\
& \geq {\cal A} \big( \{z \in \D\,;\ \Re \phi_1 (z) < \frac{\pi}{2} ( 1 - \frac{1}{h}) \} \big) \\
& = {\cal A} \big( \{z \in \D\,;\ | f (z) | < \e^{\pi/2} \e^{- \pi/2h} \} ) \\
& \geq c\, {\cal A} \big( \{z \in \D\,;\ |z| < \e^{\pi/2} \e^{- \pi/2h} \} ) = c'\, \e^\pi \e^{- \pi/h} \,,
\end{align*}
for $h > 0$ small enough.\par
\smallskip

The proof of Theorem~\ref{symbol} is completed. \qed
\bigskip

\noindent{\bf Remarks.}\par
1) In Theorem~\ref{symbol}, we have in particular $\rho_\phi (h) = O\, (h^\alpha)$ ({\it i.e.} $m_\phi$ is an 
$\alpha$-Carleson measure) for every $\alpha \geq 1$; hence (\cite{JFA}, Proposition~3.2), the composition 
operator $C_\phi \colon H^2 \to H^2$ is in all the Schatten classes $S_p (H^2)$, $p > 0$.\par\smallskip

2) In the opposite direction of Theorem~\ref{theo Hardy-Bergman}, it would be interesting to have, for every Orlicz 
function $\Psi$, a composition operator which is compact on ${\mathfrak B}^\Psi$ but not compact on $H^\Psi$. 
This is the case for $\Psi (x) = \e^{x^2}  - 1$  (\cite{CompOrli}, Theorem~5.8). Theorem~3.1 of \cite{LLQR-Rev} 
could give such an example. It could also give an example where the condition of \ref{coro compacite} is not 
sufficient to have compactness.\par
However, we may remark that the compactness of $C_\phi \colon H^\Psi \to H^\Psi$ implies the compactness of 
$C_\phi \colon {\mathfrak B}^{\Psi^2} \to {\mathfrak B}^{\Psi^2}$ since, if $\tilde \Psi (x) = [\Psi (x)]^2$, then 
$\tilde \Psi^{-1} (t^2) = \Psi^{-1} (t)$, so 
$\tilde \Psi^{-1} (1/h^2) / \tilde \Psi^{-1} \big(1/ \nu_{\phi, 2} (h) \big) \leq 
\Psi^{-1} (1/h) / \Psi^{-1} \big(1/ \nu_\phi (h) \big)$, since $\nu_{\phi, 2} (h) \leq [\nu_\phi (h)]^2$, where 
$\nu_\phi (h) = \sup_{|w| \geq 1 - h} N_\phi (w)$.\par



\bigskip

\vbox{\noindent{\it 
{\rm Pascal Lef\`evre}, Univ Lille Nord de France F-59\kern 1mm 000 LILLE, FRANCE\\
UArtois, Laboratoire de Math\'ematiques de Lens EA~2462, \\
F\'ed\'eration CNRS Nord-Pas-de-Calais FR~2956, \\
F-62\kern 1mm 300 LENS, FRANCE \\ 
pascal.lefevre@euler.univ-artois.fr 
\smallskip

\noindent
{\rm Daniel Li}, Univ Lille Nord de France F-59\kern 1mm 000 LILLE, FRANCE\\
UArtois, Laboratoire de Math\'ematiques de Lens EA~2462, \\
F\'ed\'eration CNRS Nord-Pas-de-Calais FR~2956, \\
Facult\'e des Sciences Jean Perrin,\\
Rue Jean Souvraz, S.P.\kern 1mm 18, \\
F-62\kern 1mm 300 LENS, FRANCE \\ 
daniel.li@euler.univ-artois.fr
\smallskip

\noindent
{\rm Herv\'e Queff\'elec}, Univ Lille Nord de France F-59\kern 1mm 000 LILLE, FRANCE\\
USTL, Laboratoire Paul Painlev\'e U.M.R. CNRS 8524, \\
F-59\kern 1mm 655 VILLENEUVE D'ASCQ Cedex, FRANCE \\ 
queff@math.univ-lille1.fr
\smallskip

\noindent
{\rm Luis Rodr{\'\i}guez-Piazza}, Universidad de Sevilla, \\
Facultad de Matem\'aticas, Departamento de An\'alisis Matem\'atico,\\ 
Apartado de Correos 1160,\\
41\kern 1mm 080 SEVILLA, SPAIN \\ 
piazza@us.es\par}
}

\end{document}